\documentclass[11pt]{article}
\usepackage{amsfonts}
\usepackage{latexsym, amssymb, amsmath, amscd, amsfonts, epsfig, graphicx, colordvi}
\usepackage{ifpdf}
\usepackage{graphicx}
\usepackage{xcolor}
\usepackage[dvipsnames]{pstricks}
\newtheorem{thm}{Theorem}[section]

\newtheorem{pro}[thm]{Proposition}

\newtheorem{defi}[thm]{Definition}
\newtheorem{lem}[thm]{Lemma}

\newtheorem{exam}[thm]{\bf Example}
\def\pf{\noindent{\it Proof.} }
\def\qed{\nopagebreak\hfill{\rule{4pt}{7pt}}
	\medbreak}

\def\qed{\nopagebreak\hfill{\rule{4pt}{7pt}}
	\medbreak}

\makeatletter
\def\ExtendSymbol#1#2#3#4#5{\ext@arrow 0099{\arrowfill@#1#2#3}{#4}{#5}}
\makeatother

\title{Lattice paths and Rogers--Ramanujan--Gordon type overpartitions}
\author {Diane Y. H.
	Shi}
\date{	\vspace{15pt} School of Mathematics, \\ Tianjin University, Tianjin 300072, P.R. China
	\vskip 0.2 cm
	Email: shiyahui@tju.edu.cn}

\begin{document}
	
	\maketitle
	\noindent {\bf Abstract.}	
In this paper, we establish a connection between Rogers-Ramanujan-Gordon type overpartitions to lattice paths with four kinds of unitary steps. By establishing the bijective relationship between overpartitions and lattice paths, we demonstrate that the theorems provided by Chen, Sang and Shi can be formulated in the form of lattice paths.  Subsequently, inspired by Andrews' work on parity in partition identities and  its related implications, we  impose parity constraints on lattice paths and  present some novel discoveries. By leveraging the parity outcomes within lattice paths, we also revisit overpartitions to derive results pertaining to overpartitions with parity considerations.

	\noindent {\bf Keywords:} Rogers--Ramanujan--Gordon Theorem,
	Andrews--Gordon identity, overpartition, Gordon marking, lattice path, parity
	
	\noindent {\bf AMS Subject Classification:} 05A17, 11P84	
	
	\section{Introduction}
	
The first elegant generalization of the Rogers-Ramanujan identities was given by Gordon \cite{gor61} in 1961 using a combinatorial form based on partitions. The partitions defined by Gordon that satisfy certain difference conditions are referred to as Rogers-Ramanujan-Gordon type partitions. Subsequently, Andrews \cite{and66} provided a generalization in the form of a $q$-series identity, which can also be viewed as the generating function identity of Gordon's partition theorem. We refer to it as the Andrews-Gordon identity.

Both Gordon's and Andrews' generalizations concern odd moduli. Bressoud \cite{bre79,Bre80} provided generalizations for even moduli in both combinatorial and analytic forms, which are also theorems of Rogers-Ramanujan-Gordon type partitions and Andrews-Gordon type identities.

Both Andrews-Gordon type identities can be viewed as generating functions for partitions, but this type of identity also has other combinatorial interpretations, one of which is lattice paths. In 1989, Bressoud \cite{bre87} provided a series of results on lattice paths, which included lattice path forms of the results related to Gordon, Andrews, and Bressoud.

In his 2010 paper, Andrews \cite{And10} studied the problem of adding parity restrictions to Rogers-Ramanujan-Gordon type partitions and provided results for certain parameter cases. Results for other parameter cases were subsequently given by Kur\c{s}ung\"{o}z\cite{kur09,kur10}, Kim and  Yee \cite{kim13}, while the lattice path forms of these results were presented by Hao and Shi in 2024 \cite{Shi23}.

In \cite{chen13a,chen13b,sang15}, Chen, Sang, and Shi derived the following overpartition analogues of the results by Gordon, Andrews, and Bressoud, which is a generalization of Lovejoy's results\cite{lov03}.
\begin{thm}\label{thmlast1} For $k\geq a\geq 1$, let $\overline{B}_{k,a}(n)$ denote the number of overpartitions of $n$ of the form $\lambda_1+\lambda_2+\cdots+\lambda_s$,
	such that $1$ can occur as a non-overlined part at most $a-1$ times,
	and where $\lambda_j-\lambda_{j+k-1}\geq1$ if $\lambda_j$ is overlined and
	$\lambda_j-\lambda_{j+k-1}\geq2$ otherwise. For $k> a\geq 1$, let $\overline{A}_{k,a}(n)$ denote the number of
	overpartitions of $n$ whose non-overlined parts are not congruent to
	$0,\pm a$ modulo $2k$ and let $\overline{A}_{k,k}(n)$ denote the number of overpartitions of $n$ with parts not divisible by $k$.
	Then $\overline{A}_{k,a}(n)=\overline{B}_{k,a}(n)$.
	
	And the generating function form  as the following identity:
	\begin{align}\label{AB}
		&\sum_{N_1\geq\cdots\geq N_{k-1}\geq0}
		\frac{q^{\frac{(N_1+1)N_1}{2}+N_2^2+\cdots+N_{k-1}^2+N_{a+1}+\cdots+N_{k-1}}
			(-q)_{N_1-1}(1+q^{N_a})}{(q)_{N_1-N_2}\cdots(q)_{N_{k-2}-N_{k-1}}(q)_{N_{k-1}}}\nonumber
		\\[6pt] &  \qquad \qquad =\frac{(-q)_\infty(q^a,q^{2k-a},q^{2k};q^{2k})_\infty}{(q)_\infty}.
	\end{align}
\end{thm}

\begin{thm}\label{main}For $k-1\geq a\geq 1$, let $\overline{D}_{k,a}(n)$ denote the number of overpartitions of $n$ of the form $\lambda=\lambda_1+\lambda_2+\cdots+\lambda_s$,
	such that \begin{itemize}
		\item[(i)] $f_1(\lambda)\leq a-1$, \item[(ii)] $f_l(\lambda)+f_{\overline{l}}(\lambda)+f_{l+1}(\lambda)\leq k-1$, and
		\item[(iii)] if $f_l(\lambda)+f_{\overline{l}}(\lambda)+f_{l+1}(\lambda)=k-1$, then \[lf_l(\lambda)+lf_{\overline{l}}(\lambda)+(l+1)f_{l+1}(\lambda)\equiv V_{\lambda}(l)+a-1 (\rm{mod}\; 2).\]\end{itemize} Let $\overline{C}_{k,a}(n)$ denote the number of
	overpartitions of $n$ whose non-overlined parts are not congruent to
	$0,\pm a$ modulo $2k-1$. Then $\overline{C}_{k,a}(n)=\overline{D}_{k,a}(n)$.
	
	The generating function form is as follows.	
	For $k-1 \geq a\geq 1$,
	\begin{align}&\sum_{N_1\geq N_2\geq\cdots\geq
			N_{k-1}\geq0}\frac{q^{\frac{(N_1+1)N_1}{2}+N_2^2+\cdots+N_{k-1}^2+N_{i+1}+\cdots+N_{k-1}}
			(-q)_{N_1-1}(1+q^{N_a})}{(q)_{N_1-N_2}\cdots(q)_{N_{k-2}-N_{k-1}}(q^2;q^2)_{N_{k-1}}}\nonumber\\ \label{eqD}&=\frac{(-q)_\infty(q^a,q^{2k-1-a},q^{2k-1};q^{2k-1})_\infty}{(q)_\infty}.\end{align}
\end{thm}
The identity \eqref{AB} is derived by Chen, Sang and Shi \cite{chen13a} and the identity \eqref{eqD} is derived by Sang and Shi \cite{sang15}.
The left hand side of \eqref{AB} and \eqref{eqD} are in infinite summation forms which are derived from the Gordon marking for overpartitions.

In this paper, we  focus on the combinatorial object -- 
lattice paths. We shall establish the relationship between overpartitions and lattice paths with four types of unit steps. Then, we will provide lattice path interpretations for the two aforementioned identities.  Additionally, we present some new findings on lattice paths and Rogers--Ramanujan--Gordon type overpartitions that incorporate parity considerations.

We first introduce some definitions about overpartition and lattice path.
	A partition $\lambda$ of a positive  integer $n$ is a non-increasing sequence of positive integers (called the parts of $\lambda$) $\lambda_1\geq \cdots\geq \lambda_s>0$ such that $n=|\lambda|=\lambda_1+\cdots+\lambda_s$. An overpartition $\lambda$ of a positive integer $n$ is also a non-increasing sequences of positive integers $\lambda_1\geq \cdots\geq \lambda_s>0$ such that $n=|\lambda|=\lambda_1+\cdots+\lambda_s$ and the first occurrence  of  each integer may be overlined.  The partition or overpartition of integer zero is $\emptyset$, so the number of partitions and overpartitions of integer $0$ are both $1$. The number of parts $s=s(\lambda)$ of $\lambda$ is the length of $\lambda$. For example, $(\overline{7},6,6,\overline{5},2,\overline{1})$ is an overpartition of $27$  with $s=6$.
	Given a partition or an overpartition $\lambda$, let $f_l(\lambda) (f_{\overline{l}}(\lambda))$ denote the number of occurrences of 
	non-overlined (overlined) $l$ in $\lambda$. Let $V_{\lambda}(l)$ denote the number of overlined parts in $\lambda$ less than or
	equal to $l$.
	
	We use the common notation in $q$-series. Let
\[(a)_\infty=(a;q)_\infty=\prod_{i=0}^{\infty}(1-aq^i),\]
and
\[(a)_n=(a;q)_n=(1-a)(1-aq)\cdots(1-aq^{n-1}).\]
We also write
\[(a_1,\ldots,a_k;q)_\infty=(a_1;q)_\infty\cdots(a_k;q)_\infty.\]

In this study, we examine lattice paths situated within the first quadrant, which start on the $y$-axis, end on the $x$-axis, 	and incorporate four distinct types of unit steps:
\begin{itemize}
	\item 	North–East NE: $(x,y) \rightarrow  (x + 1,y + 1)$,
	
\item	South–East SE: $(x,y) \rightarrow (x+1,y-1)$,
	
\item	 North N: $(x,y-1)\rightarrow (x,y)$,
\item	East E: $(x,0)\rightarrow (x+1,0)$.
\end{itemize}	
In the lattice paths, an East step can only appear at height 0.   

\begin{figure}[h]
\begin{center}
	\begin{picture}(150,50)
\thicklines
\put(10,20){\line(0,1){25}}\put(10,45){\line(1,-1){25}}
\put(95,20){\line(1,1){25}}\put(120,45){\line(1,-1){25}}
\put(0,0){NSE peak} \put(95,0){NESE peak}
\end{picture}
\end{center}\caption{The two kinds of peaks.}\label{fig:figureP}	\end{figure}
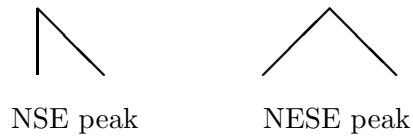

Our lattice paths with four types of unit steps differ from those of  Corteel and  Mallet \cite{cor07}; however, we will adopt some of the parameter definitions from their article.
The height of a vertex corresponds to the $y$ coordinate, and the weight of a vertex is its $x$ coordinate. A peak $(x,y)$ is a vertex that is preceded by a North-East step (in which case it is called an NESE peak) or a North step (in which case it is called an NSE peak) and is followed by a South-East step (Figure \ref{fig:figureP}). 
These paths are either empty or end with a South-East step. If the leftmost peak is an NSE peak, it must start with an NE step, except when it starts at $(0,0)$ with an E step.
The major index of a path is the sum of the weights of its peaks. Let $k$ and $a$ be positive integers with $a \leq k$. We say that a path satisfies the special $(k,a)$-conditions if it starts at height $k-a$ and its height is less than $k$.

We present  the main results of this paper.
\begin{thm}\label{EB}
For $k \geq a\geq 1$, let $\overline{E}_{k,a}(n)$ be the number of paths with the four kinds steps introduced above whose  major index is $n$ satisfying the special $(k,a)$-conditions. Then we have 

\[\overline{E}_{k,a}(n)=\overline{B}_{k,a}(n),\] and the generating function of $\overline{E}_{k,a}(n)$ is \eqref{AB}.
\end{thm}

We revisited the definition of the relative height of a peak, which was proposed by Bressound in \cite{bre87} and elaborated by Corteel and Mallet in \cite{cor07}. Bressoud's definition is in the context of paths with three basic steps, whereas the paths defined by Corteel and Mallet in \cite{cor07}, which include four basic steps, differ from the paths defined in this paper that also incorporate four basic steps. However, the definition of relative height remains similar, with only minor differences.

\begin{defi}
The relative height of a peak $(x,y)$ is the largest integer $h$ for which we can find
		two vertices on the path, $(x',y'-h)$ and $(x'',y-h)$, such that $x'\leq x<x''$ and such that between
		these two vertices there are no peaks of height $>y$ and every peak of height $y$ has  abscissa $\geq x$.
\end{defi}
For example, the four peaks in the following lattice path have relative height of $1,\ 1,\ 3$ and $1$ respectively.
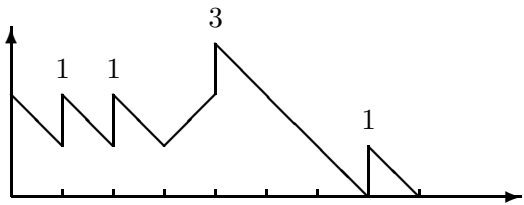
\begin{center}
	\begin{picture}(180,80)
		\thicklines
		\setlength{\unitlength}{0.45mm}
		\multiput(25,10)(15,0){8}{\line(0,1){2}}
		\put(10,10){\vector(1,0){150}}
		\put(10,10){\vector(0,1){50}}
\put(10,40){\line(1,-1){15}}\put(25,25){\line(0,1){15}}\put(25,40){\line(1,-1){15}}\put(40,25){\line(0,1){15}}\put(40,40){\line(1,-1){15}}\put(55,25){\line(1,1){15}}\put(70,40){\line(0,1){15}}\put(70,55){\line(1,-1){45}}\put(115,10){\line(0,1){15}}\put(115,25){\line(1,-1){15}}\put(23,45){$1$}\put(38,45){$1$}\put(68,60){$3$}\put(113,30){$1$}
\end{picture}\end{center}
Let $W(x)$ denote the number of NSE peaks whose weight is less than or equal to $x$. In the above lattice path, $W(1)=1,\ W(2)=W(3)=2,\ W(4)=W(5)=W(6)=3$ and $W(7)=W(8)=4.$
\begin{thm}\label{FD}
For $k \geq a\geq 1$, let $\overline{F}_{k,a}(n)$ denote the number of lattice paths with major index $n$ that satisfy the special $(k,a)$-conditions, using the four kinds of unitary steps introduced above, and for peaks with relative height $k-1$, the weight $x$ satisfies \begin{equation}\label{equiv}x \equiv W(x)+a-1 \pmod{2}.\end{equation} Then, we have \[\overline{F}_{k,a}(n)=\overline{D}_{k,a}(n),\] and the generating function of $\overline{F}_{k,a}(n)$ is given by \eqref{eqD}.
\end{thm}

The lower even parity index in partition was first introduced by Andrews \cite{And10}, and then Kur\c{s}ung\"{o}z \cite{kur10} used it to provide a combinatorial interpretation of a formula proposed by Andrews in \cite{And10} as an open problem. We have defined the  parity of a peak and the lower even parity index in lattice path with three unitary steps in \cite{Shi23}. In this parper, we also introduce these concepts in lattice paths with four kinds of unitary steps. 

\begin{defi}
For a peak $(x,y)$ with relative height $r$, we define the parity as the opposite parity of  $x-r-W(x)$.
\end{defi}
	Next, we will provide some results about lattice paths with parity restrictions on peaks.
\begin{thm}\label{thmH}
For $k \geq a \geq 1$, let $H_{k,a}(n)$ denote the number of lattice paths that satisfy the special $(k,a)$-conditions and have a major index of $n$ with all peaks having even parity. Then we can obtain the generating function for $H_{k,a}(n)$ as follows.
\small{\begin{align}
		&\sum_{n\geq 0}H_{k,a}(n)q^n\nonumber\\&=\sum_{n_1,\ldots,n_{k-1}\geq 0}	\frac{q^{N_1(N_1+1)/2+N_2^2+\cdots+N_{k-1}^2+2N_{a+1}+2N_{a+3}+\cdots+2N_{k-1}}(-q)_{N_1}}{(q^2;q^2)_{n_1}\cdots(q^2;q^2)_{n_{k-2}}(q^2;q^2)_{n_{k-1}}}.
\end{align}}
Throughout the entire paper, we assume that $N_k=0$ and  $n_i= N_i-N_{i+1}$, for $1\leq i\leq k-1$.
\end{thm}

\begin{thm}\label{thmJ}
For $k \geq a\geq 1$, let $J_{k,a}(n)$ denote the number of lattice paths that satisfy the special $(k,a)$-conditions and have a major index of $n$ with all peaks having odd parity. Then, we have  the generating function of $J_{k,a}(n)$ is 
\small{\begin{align}
		&\sum_{n\geq 0}J_{k,a}(n)q^n\nonumber\\&\label{eqJ}=\sum_{n_1,\ldots,n_{k-1}\geq 0}	\frac{q^{N_1(N_1+1)/2+N_2^2+\cdots+N_{k-1}^2+2N_{a+2}+2N_{a+4}+\cdots+2N_{k-2}}(-q)_{N_1-1}(q^{N_1}+q^{2N_a})}{(q^2;q^2)_{n_1}\cdots(q^2;q^2)_{n_{k-2}}(q^2;q^2)_{n_{k-1}}}.
\end{align}}
\end{thm}

In \cite{And10}, one of the open problems Andrews listed is to investigate the three-parameter function 
\begin{equation}\label{eqARRG}
	\sum_{N_1\geq N_2\geq \cdots \geq N_{k-1}\geq 0}\frac{q^{N_1^2+N_2^2+\cdots+N_{k-1}^2}x^{N_1+N_2+\cdots+N_{k-1}}(-yq)_{n_
			1}(-yq)_{n_2}\ldots(-yq)_{n_{k-1}}}
	{(q^2;q^2)_{n_1}\ldots(q^2;q^2)_{n_{k-2}}(q^2;q^2)_{n_{k-1}}}.
\end{equation} 

Kurşungöz explained the formula \eqref{eqARRG} using clusters in Gordon marking in his paper \cite{kur10}. Hao and Shi interpreted it through lattice paths in \cite{Shi23}.

We provide a three-parameter function analogous to the formula \eqref{eqARRG}, which can be interpreted by overpartitions and lattice paths.

\begin{thm}\label{thmGT}For $k\geq 1$, let $G_{k,k}(l,m,n)$ dennote the number of lattice paths enumerated by $\overline{E}_{k,k}(n)$ whose sum of the relative height of peaks is $m$ and whose full lower even peak parity index is $l$. Let $T_{k,k}(l,m,n)$ denote the number of overpartitions enumerated by $\overline{B}_{k,k}(n)$ which have $m$ parts and whose full lower even cluster parity index is $l$.  Consequently, it is established that
	\begin{equation}\label{GT}G_{k,k}(l,m,n)=T_{k,k}(l,m,n),\end{equation}and the generative function corresponding to $G_{k,k}(l,m,n)$ and $T_{k,k}(l,m,n)$ can be stated as: 
	\begin{align}
		&\nonumber\sum_{l,m,n\geq 0}G_{k,k}(l,m,n)y^lx^mq^n=\sum_{l,m,n\geq 0}T_{k,k}(l,m,n)y^lx^mq^n\\=&\sum_{N_1\geq\cdots\geq N_{k-1}\geq0}
		\frac{q^{\frac{(N_1+1)N_1}{2}+N_2^2+\cdots+N_{k-1}^2}
			(-q)_{N_1-1}(1+q^{N_k})x^{N_1+\cdots+N_{k-1}}}{(q^2;q^2)_{n_1}\cdots(q^2;q^2)_{n_{k-2}}(q^2;q^2)_{n_{k-1}}}\nonumber\\&\label{equ2}\qquad \times (-yq)_{n_1}(-yq)_{n_2}\cdots(-yq)_{n_{k-1}}.
\end{align}\end{thm}
It is worth noting that when $N_k = 0$, we replace 1 with $q^{N_k}$ to maintain consistency with the previous equation format.

The definitions of full lower even peak parity index and full lower even cluster parity index will be stated in Section 2.

 In Section 2, we introduce the background of Gordon marking for overpartitions and also present some new definitions related to overpartitions and lattice paths containing four types of unit steps. Following that, we provide some new results about overpartitions, which can be seen as analogues of Andrews' parity results for overpartitions. And we recall some maps in \cite{chen13a}.  In Section 3, we present the bijection between the Gordon marking of overpartitions and the lattice paths.
In Section 4--6, we give the proofs of Theorem \ref{EB}, \ref{FD}, \ref{thmH}, \ref{thmJ}, \ref{thmO}, \ref{thmoverO} and \ref{thmGT}.

\section{Background and preliminaries}

In this section, we shall introduce some background about the Gordon marking of overpartition  \cite{chen13a,chen13b,sang15}. Then, we list some results which can be seen as the overpartition analogues of Andrews' parity results.  Additionally, we also provide the definitions of the parity index in overpartition and lattice path, which are used in Theorem \ref{thmGT}.
In the final part of this section, we present some operations on lattice paths, which will be used in the bijective constructions in subsequent chapters. 

Let  $\mathcal{\overline{B}}_{k,a}(n)$
denote the set of overpartitions enumerated by $\overline{B}_{k,a}(n)$. Let $\mathcal{U}_{k,a}(n)$ (resp. $\mathcal{V}_{k,a}(n)$)
denote the subset of $\mathcal{\overline{B}}_{k,a}(n)$ in which  the overpartitions with the smallest parts have an
overlined part (resp. have no overlined part), and let $U_{k,a}(n)$ (resp. $V_{k,a}(n)$)  denote the number of overpartitions in $\mathcal{U}_{k,a}(n)$ (resp. $\mathcal{V}_{k,a}(n)$). 

We revisit the Gordon marking for overpartition initially presented by Chen, Sang and Shi \cite{chen13a}.

\begin{defi}
	The Gordon marking of an overpartition $\lambda$ is an assignment of
	positive integers (marks) to parts of $\lambda$. We assign the marks to parts
	in the following order
	\begin{equation}\label{order}
		\overline{1} < 1 < \overline{2} < 2< \cdots
	\end{equation}
	such that the marks are as small as possible subject to the following conditions.
	If $\overline{j+1}$ is not a part of $\lambda$,
	then all the parts $j$, $\overline{j}$, and $j+1$ are assigned  different integers.
	If  $\lambda$ contains an overlined part $\overline{j+1}$,
	then the smallest mark
	assigned to a part $j$ or $\overline{j}$ can be used as the mark of  $j+1$ or $\overline{j+1}$.
\end{defi}

For example, given an overpartition
\[\lambda= (16,13,12,12,11,\overline{10},\overline{8},8,8,7,\overline{6},6,5,5,\overline{4},2,2,\overline{1}).\]
The
Gordon marking of $\lambda$ is
\[ (\overline{1}_1,2_2,2_3,\overline{4}_1,5_2,5_3,
\overline{6}_1,6_2,7_3,\overline{8}_1,8_2,8_3,\overline{10}_1,11_2,12_1,12_3,13_2,16_1),
\]
where the subscripts are the marks. The Gordon marking of $\lambda$ can
also be illustrated as
\begin{equation}\nonumber\label{lambda}
	\lambda=\setcounter{MaxMatrixCols}{16}\begin{bmatrix}
		\ &2&\ &\ &5&\ &7&8&\ &\ &\ &12&\ &\ &\ &\\[3pt]
		\ &2&\ &\ &5&6&\ &8&\ &\ &11&\ &13&\ &\ &\\[3pt]
		\overline{1}&\ &\ &\overline{4}&\ &\overline{6}&\ &\overline{8}&\ &\overline{10}&\ &12&\ &\ &\ &16
	\end{bmatrix}\begin{matrix}
		3\\[3pt]2\\[3pt]1
	\end{matrix}\;,
\end{equation}
where
the parts in the third  row are marked by $1$, the parts in the second row
are marked by $2$, and the parts in the first row are marked by $3$.

It can be observed that the Gordon marking for any overpartition is unique. It can also be verified that all the overlined parts are marked with $1$.

Let us revisit some definitions and mappings defined in \cite{chen13a}. Let $\overline{\mathcal{B}}_{N_1, N_2, \ldots, N_{k-1}; a}(n)$ (respectively, $\mathcal{U}_{N_1, N_2, \ldots, N_{k-1}; a}(n)$ and $\mathcal{V}_{N_1, N_2, \ldots, N_{k-1}; a}(n)$) denote the subset of $\overline{\mathcal{B}}_{k, a}(n)$ (respectively, $\mathcal{U}_{k, a}(n)$ and $\mathcal{V}_{k, a}(n)$) in which the number of parts marked with $i$ is exactly $N_i$, for $i = 1, \ldots, k-1$. Meanwhile, $\overline{B}_{N_1, N_2, \ldots, N_{k-1}; a}(n)$ (respectively, $U_{N_1, N_2, \ldots, N_{k-1}; a}(n)$ and $V_{N_1, N_2, \ldots, N_{k-1}; a}(n)$) represents the number of overpartitions in $\overline{\mathcal{B}}_{N_1, N_2, \ldots, N_{k-1}; a}(n)$ (respectively, $\mathcal{U}_{N_1, N_2, \ldots, N_{k-1}; a}(n)$ and $\mathcal{V}_{N_1, N_2, \ldots, N_{k-1}; a}(n)$). Let $\mathcal{P}_{N_1, N_2, \ldots, N_{k-1}; a}(n)$ denote the set of overpartitions in $\mathcal{U}_{N_1, N_2, \ldots, N_{k-1}; a}(n)$ where all parts marked with $1$ are overlined.
Set
 \begin{align}
 \nonumber 	&\mathcal{U}_{N_1,N_2,\ldots,N_{k-1};a}=\bigcup_{n\geq 0}\mathcal{U}_{N_1,N_2,\ldots,N_{k-1};a}(n),\\
 \nonumber	&\mathcal{P}_{N_1,N_2,\ldots,N_{k-1};a}=\bigcup_{n\geq 0}\mathcal{P}_{N_1,N_2,\ldots,N_{k-1};a}(n).
 \end{align}
Let $m$ denote the number of parts in overpartitions in $\overline{\mathcal{B}}_{N_1,N_2,\ldots,N_{k-1};a}(n)$, then $m=N_1+N_2+\cdots+N_{k-1}$.
For an overpartition $\lambda\in \mathcal{U}_{N_1,N_2,\ldots,N_{k-1};a}$, according to the bijection $\varphi$ \cite[Section 4]{chen13a},  $\lambda$ corresponds to a distinct partition $\beta$ with at most $N_1-1$  parts and an overpartition $\alpha\in \mathcal{P}_{N_1,N_2,\ldots,N_{k-1};a}$, where $|\lambda|=|\alpha|+|\beta|$. 
Since $k-1$ is the largest mark in $\lambda$, we define the set $\mathcal{Q}_{N_1,N_2,\ldots,N_{k-1};a}(n)$ as the set of overpartitions $\lambda$ in $\mathcal{P}_{N_1,N_2,\ldots,N_{k-1};a}(n)$  that satisfy the following relation,
\begin{equation}\label{fte}
	f_t(\lambda)+f_{\overline{t}}(\lambda)+f_{t+1}(\lambda) =  k-1,
\end{equation}
for any positive integer $t$ that is smaller than the greatest $(k-1)$-marked part. Set
\[\mathcal{Q}_{N_1,N_2,\ldots,N_{k-1};a}=\bigcup_{n\geq 0}\mathcal{Q}_{N_1,N_2,\ldots,N_{k-1};a}(n).\]

In \cite[Section 5]{chen13a}, by using the second reduction and dilation operations, the authors construct   a bijection $\psi$ between
$\mathcal{P}_{N_1,N_2,\ldots,N_{k-1};a}$ and $\mathcal{Q}_{N_1,N_2,\ldots,N_{k-1};a}\times R_{n_{k-1}}$,
where $R_{n_{k-1}}$ denotes the set of ordinary partitions with  $n_{k-1}$ nonnegative parts. 

We will further utilize the bijective function $\chi$  \cite[Section 6]{chen13a} which acts  between $\mathcal{Q}_{N_1,\ldots,N_{k-1};a}(n)$
and $\mathcal{Q}_{N_1-1,\ldots,N_{k-1}-1;a}(n-N_1-2N_2-\cdots-2N_{k-1}+a-1)$. 

For convenience, we introduce the definition of a cluster in overpartitions, which is identical to the definition of a cluster in partitions \cite{kur10}.
\begin{defi}	
	In the overpartition $\lambda=\lambda_1+\cdots+\lambda_m$, an $r$-cluster is a sub-partition $\lambda_{i_1}\leq \cdots \leq \lambda_{i_r}$ such that $\lambda_{i_j}$ is $j$-marked for $j=1, \ldots, r$; $\lambda_{i_{j+1}}-\lambda_{i_j}\leq 1$ for $j=1, \ldots, r-1$, and there are no $(r+1)$-marked parts equal to $\lambda_{i_r}$ or $\lambda_{i_{r}}+1$.
\end{defi}

However, the cluster decomposition is not unique in the Gordon marking of overpartitions. We provide an example to illustrate that the cluster decomposition in overpartition Gordon marking may not be unique.
\begin{exam}
\begin{equation}
\lambda=	\setcounter{MaxMatrixCols}{15}\begin{bmatrix}\nonumber
		\ &2&\ &\ &5&\ &7&8&\ &\ &\ &\mathbf{12}&\ &\ &\ \\[3pt]
		\ &2&\ &\ &5&6&\ &8&\ &\ &\mathbf{11}&\mathbf{12} &\ &\ &\ \\[3pt]
		\overline{1}&\ &\ &4&\ &\overline{6}&\ &\overline{8}&\ &\mathbf{\overline{10}}&\ &\mathbf{\overline{12}}&\ &\ &15
	\end{bmatrix}\begin{matrix}
		3\\[3pt]2\\[3pt]1
	\end{matrix}\;.
\end{equation}

We can see that the $3$-marked  $12$ can be in a $3$-cluster $(\overline{10},11,12)$ or $(\overline{12},12,12)$. This is because there is a $2$-marked $11$ and a $2$-marked $12$. 
\end{exam}

If we decompose the overpartition from left to right, making the clusters as long as possible and ensuring that each part belongs to only one cluster, then we can obtain a decomposition into $N_1$ clusters, referred to as the left clusters. If we decompose the Gordon marking from right to left, making the clusters as long as possible, we may obtain a different decomposition, referred to as the right clusters.

In nearly all the instances within this paper, we employ the left cluster decomposition, hence the term ``left" is omitted. Moreover, a left $r$-cluster with the smallest part being $l$ is denoted as $C_r(l)$. It can be directly deduced that an overpartition in $\overline{\mathcal{B}}_{N_1, N_2, \ldots, N_{k-1}; a}(n)$ possesses $n_i$ $i$-clusters. Notably, the subtraction of parts in $\chi$ essentially represents a $(k-1)$-cluster with the smallest $(k-1)$-marked part.
 
The left cluster decomposition of the above example is $3$-cluster $(\overline{1},2,2)$,  $(4,5,5)$,  $(\overline{6},6,7)$,  $(\overline{8},8,8)$,  $(\overline{10},11,12)$, denoted as $C_3(1)$, $C_3(4)$, $C_3(6)$, $C_3(8)$ and $C_3(10)$ respectively. $2$-cluster $(\overline{12},12)$ and $1$-cluster $(15)$, denoted as $C_2(12)$ and $C_1(15)$.  The right cluster decomposition includs $3$-cluster $(\overline{1},2,2)$,  $(4,5,5)$,  $(\overline{6},6,7)$,  $(\overline{8},8,8)$, $2$-cluster $(\overline{10},11)$, $3$-cluster $(\overline{12},12,12)$ and $1$-cluster $(15)$.

\begin{defi}
	For an $r$-cluster $C_r(l)$ in the Gordon marking of an overpartition, its parity is defined as the opposite parity of the number of even parts in the cluster minus $V(l)$.
\end{defi}
It is worth mentioning that is  
\begin{equation}\label{clusterparity}
	\text{the sum of the parts in cluster} -r\equiv \text{the number of even parts in  cluster}\pmod{2}.\end{equation}

We provide the overpartition analogues of Andrews' parity results.
\begin{thm}\label{thmO}
	Let $O_{k,a}(n)$ denote the number of overpartitions enumerated by $\overline{B}_{k,a}(n)$ such that the number of even parts in each cluster minus $V(l)$ is even, i.e., the  clusters all have  odd parity in the overpartitions, where $l$ is the $1$-marked part in the cluster. Then for $1\leq a < k$, 
	\[O_{k,a}(n)=J_{k,a}(n).\]
	we have  the generating function is
\small{\begin{align}
\sum_{n\geq 0}O_{k,a}(n)q^n=\sum_{n_1,\ldots,n_{k-1}\geq 0}	\frac{q^{N_1(N_1+1)/2+N_2^2+\cdots+N_{k-1}^2+2N_{a+2}+\cdots}(q^{2N_a}+q^{N_1})(-q)_{N_1-1}}{(q^2;q^2)_{n_1}\cdots(q^2;q^2)_{n_{k-2}}(q^2;q^2)_{n_{k-1}}}.
	\end{align}}	  
\end{thm}
\begin{thm}\label{thmoverO}
	Let $\overline{O}_{k,a}(n)$ denote the number of overpartitions enumberated by $\overline{B}_{k,a}(n)$ such that the number of odd parts in each cluster minus $V(l)$ is even, where $l$ is the $1$-marked part in the cluster. Then, if $a$ is odd, the generating function is
\small{\begin{align}
&\sum_{n\geq 0}\overline{O}_{k,a}(n)q^n\nonumber\\&\label{overOo}=\sum_{n_1,\ldots,n_{k-1}\geq 0}	\frac{q^{N_1(N_1+1)/2+N_2^2+\cdots+N_{k-1}^2+N_{a+1}+\cdots+N_{k-1}}(q^{N_a+n_1+n_3+\cdots+n_{a-2}}+q^{n_2+n_4+\cdots+n_{a-1}})(-q)_{N_1-1}}{(q^2;q^2)_{n_1}\cdots(q^2;q^2)_{n_{k-2}}(q^2;q^2)_{n_{k-1}}}.
	\end{align}}
	If $a$ is even, the generating function is
	\small{\begin{align}
			&\sum_{n\geq 0}\overline{O}_{k,a}(n)q^n\nonumber\\&\label{overOe}=\sum_{n_1,\ldots,n_{k-1}\geq 0}	\frac{q^{\frac{N_1(N_1+1)}{2}+N_2^2+\cdots+N_{k-1}^2+N_a+N_{a+1}+\cdots+N_{k-1}}(q^{N_{a+1}+n_1+n_3+\cdots+n_{a-1}}+q^{n_2+n_4+\cdots+n_{a-2}})(-q)_{N_1-1}}{(q^2;q^2)_{n_1}\cdots(q^2;q^2)_{n_{k-2}}(q^2;q^2)_{n_{k-1}}}.
	\end{align}}
\end{thm}

Now we will provide the definition of the full lower even peak parity index in both overpartition and lattice path.

\begin{defi}
The lower even $r$-peak parity index of a lattice path $L$ denotes the number of times the parity of peaks with a relative height of $r$ changes from the leftmost to the rightmost peaks, starting with an even $ r$-peak if the path begins with an NSE peak; otherwise, it starts with an odd peak.
\end{defi}
\begin{defi}
	Given a lattice path $L$ that satisfies the special $(k,k)$-conditions, the full lower even peak parity index of $L$ is the sum of all lower even   $1$-, $2$-,$\ldots$, $(k-1)$-peaks parity indices.
\end{defi}
\begin{defi}
	The lower even $r$-cluster parity index of an overpartition is the number of times that the parity of $r$-clusters changes from the leftmost to the rightmost, beginning with an even $r$-cluster if the first cluster have an overlined part, otherwise begining with an odd $r$-cluster.
\end{defi}

\begin{defi}
	Given an overpartition $\lambda=(\lambda_1, \cdots,  \lambda_m)$ enumerated by $\overline{B}_{k,k}(m,n)$, the full lower even cluster parity index of $\lambda$ is the sum of all lower even $1$-, $2$-,$\ldots$,$(k-1)$-cluster parity indices.
\end{defi}

To provide the proof, we will introduce a mapping $\phi$, which maps the overpartition $\lambda$ to a lattice path $L$ composed of four types of unit steps, with its major index as $|\lambda|$. The number of NSE peaks corresponds to the number of parts with overlines in $\lambda$, and the number of peaks at relative height $i$ is $n_i$, which equals the number of $i$-clusters in $\lambda$. To construct the mapping $\phi$, we will utilize the mappings $\varphi$ \cite[Section 4]{chen13a}, $\psi$ \cite[Section 5]{chen13a}, and $\chi$ \cite[Section 6, proof of Theorem 6.2]{chen13a}.

The remaining part of this section will introduce operations on lattice paths. The first dilation operation on lattice paths is a mapping that changes an NSE peak into an NESE peak and can increase the major index of the path $L$ by one.

\noindent{\bf The first dilation operation on lattice path.}  

{\bf Type A}: The rightmost peak is an NSE peak, and the dilation operation  will directly change this NSE peak to an NESE peak.

{\bf Type B}:
Select the NSE peak whose nearest peak to the right is the rightmost NESE peak. The dilation operation turns this NSE peak into an NESE peak and the rightmost NESE peak into an NSE peak.

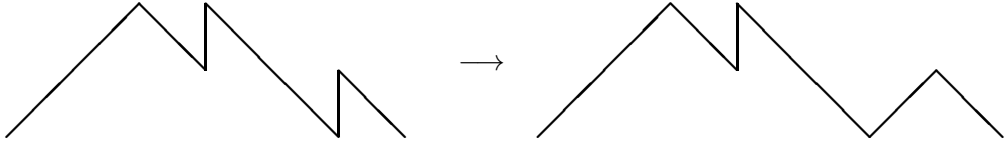
\begin{figure}[t]\centering
	\begin{picture}(400,100)
		\thicklines
		\put(0,90){Type A}	\put(10,20){\line(1,1){25}}\put(35,45){\line(1,1){25}}\put(60,70){\line(1,-1){25}}\put(85,45){\line(0,1){25}}\put(85,70){\line(1,-1){50}}\put(135,20){\line(0,1){25}}\put(135,45){\line(1,-1){25}} \put(180,45){$\longrightarrow$}\put(210,20){\line(1,1){25}}\put(235,45){\line(1,1){25}}\put(260,70){\line(1,-1){25}}\put(285,45){\line(0,1){25}}\put(285,70){\line(1,-1){50}}\put(335,20){\line(1,1){25}}\put(360,45){\line(1,-1){25}} 
	\end{picture}\caption{The first dilation operations on lattice path of Type A.} \label{figure1}\end{figure}
\begin{figure}[t]\centering
	\begin{picture}(400,100)
		\thicklines
		\put(0,90){Type B}	\put(10,20){\line(1,1){25}}\put(35,45){\line(1,1){25}}\put(60,70){\line(1,-1){25}}\put(85,45){\line(0,1){25}}\put(85,70){\line(1,-1){25}}\put(110,45){\line(1,1){25}}\put(135,70){\line(1,-1){50}} \put(180,45){$\longrightarrow$}\put(210,20){\line(1,1){25}}\put(235,45){\line(1,1){25}}\put(260,70){\line(1,-1){25}}\put(285,45){\line(1,1){25}}\put(310,70){\line(1,-1){25}}\put(335,45){\line(0,1){25}}\put(335,70){\line(1,-1){50}} 
	\end{picture}\caption{The first dilation operations on lattice path of Type B.} \label{figure2}\end{figure}
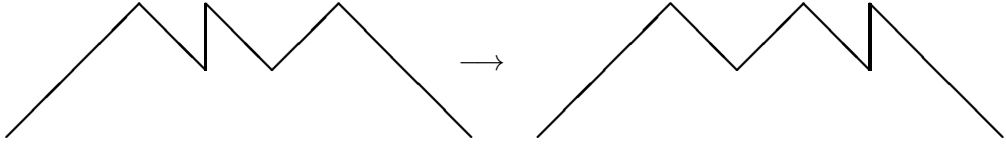
The two types dilation operations are displaied in Figure \ref{figure1} and \ref{figure2}.

The two types of the first reduction operation are the inverse of the two types of the first dilation operation, which convert an NESE peak into an NSE peak and decrease the major index by one. Since their construction is not difficult, we omit them here.

Now we will introduce the second dilation operation, which acts on an NSE peak or an NESE peak and increases the weight of the peak by one.

\noindent{\bf The second dilation operation in lattice path}

The second dilation operation applies to lattice paths whose peaks are all NSE peaks, except for the leftmost peak which can be either an NSE or an NESE peak. For a peak with a relative height of one, we increase its weight by moving the peak to the right. We show how to increase its weight by displaying them in Figure \ref{figure3} and \ref{figure3+}.

\begin{figure}[h]\centering	\begin{picture}(200,260)
		\thicklines
		\put(10,60){\line(1,-1){25}}\put(35,35){\line(0,1){25}}\put(35,60){\line(1,-1){50}}\put(60,35){\line(0,1){2}}\put(85,35){$\longrightarrow$}
		\put(110,60){\line(1,-1){50}}\put(135,35){\line(0,1){2}}
		\put(160,10){\line(0,1){25}}\put(160,35){\line(1,-1){25}}
		
		\put(15,80){\line(0,1){25}}\put(15,105){\line(1,-1){25}}\put(40,80){\line(1,0){25}}
		\put(85,90){$\longrightarrow$}\put(140,80){\line(0,1){25}}\put(140,105){\line(1,-1){25}}\put(115,80){\line(1,0){25}}	
		
		\put(10,130){\line(0,1){25}}\put(10,155){\line(1,-1){25}}\put(35,130){\line(1,1){50}}\put(60,155){\line(0,1){2}}\put(85,145){$\longrightarrow$}
		\put(110,130){\line(1,1){25}}\put(135,155){\line(0,1){25}}\put(135,180){\line(1,-1){25}}\put(160,155){\line(1,1){25}}
		
		\put(10,200){\line(0,1){25}}\put(10,225){\line(1,-1){25}}\put(35,200){\line(1,1){25}}\put(60,225){\line(0,1){25}}\put(60,250){\line(1,-1){25}}\put(85,215){$\longrightarrow$}
		\put(110,200){\line(1,1){25}}\put(135,225){\line(0,1){25}}\put(135,250){\line(1,-1){25}}\put(160,225){\line(0,1){25}}\put(160,250){\line(1,-1){25}}
	\end{picture}\caption{The rules for increaseing the weight of a NSE peak with relative height one by one} \label{figure3}\end{figure}

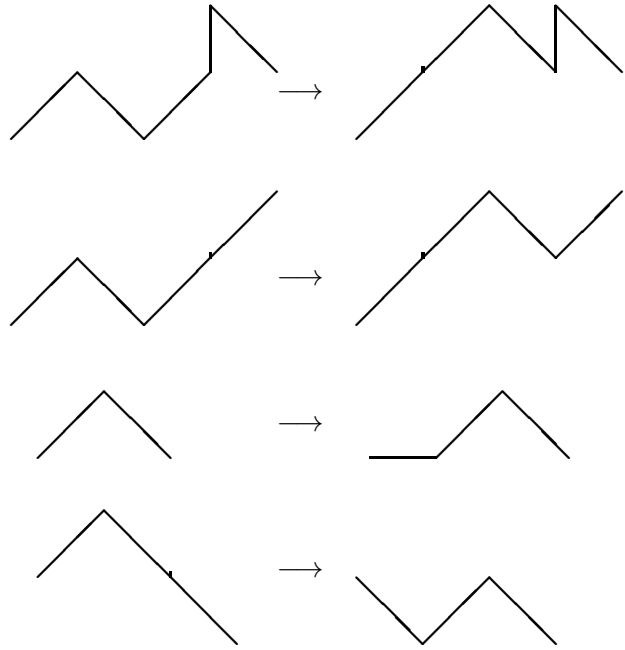
\begin{figure}[h]\centering	\begin{picture}(250,260)
		\thicklines
		\put(10,35){\line(1,1){25}}\put(35,60){\line(1,-1){50}}\put(60,35){\line(0,1){2}}\put(100,35){$\longrightarrow$}
		\put(130,35){\line(1,-1){25}}\put(155,10){\line(1,1){25}}
		\put(180,35){\line(1,-1){25}}	
		
		\put(10,80){\line(1,1){25}}\put(35,105){\line(1,-1){25}}\put(100,90){$\longrightarrow$}\put(160,80){\line(1,1){25}}\put(185,105){\line(1,-1){25}}\put(135,80){\line(1,0){25}}	
		
		\put(0,130){\line(1,1){25}}\put(25,155){\line(1,-1){25}}\put(50,130){\line(1,1){50}}\put(75,155){\line(0,1){2}}\put(100,145){$\longrightarrow$}
		\put(130,130){\line(1,1){50}}\put(155,155){\line(0,1){2}}\put(180,180){\line(1,-1){25}}\put(205,155){\line(1,1){25}}
		
		\put(0,200){\line(1,1){25}}\put(25,225){\line(1,-1){25}}\put(50,200){\line(1,1){25}}\put(75,225){\line(0,1){25}}\put(75,250){\line(1,-1){25}}\put(100,215){$\longrightarrow$}
		\put(130,200){\line(1,1){50}}\put(155,225){\line(0,1){2}}\put(180,250){\line(1,-1){25}}\put(205,225){\line(0,1){25}}\put(205,250){\line(1,-1){25}}
	\end{picture}\caption{The rules for increasing the weight of an NESE peak with relative height one by one} \label{figure3+}\end{figure}

It should be noted that when we increase the weight of a peak with relative height one, whether it is an NSE peak or an NESE peak, the peak moves to the right and may encounter the next NSE peak on the right. To be more precise, if two peaks $(x,y)$ and $(x',y')$ satisfy the condition $x'-x=1$, they will meet when the weight of $(x,y)$ is increased. If this situation occurs, the two peaks will exchange their relative heights, and subsequently, the peak on the right will become a peak with a relative height one. We then abandon the peak we were moving and proceed to move the next peak to the right.

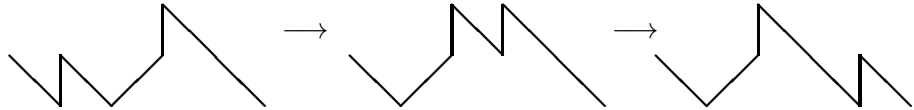
\begin{figure}[h]
	\centering\begin{picture}(340,80)
		\setlength{\unitlength}{0.45mm}\thicklines	\put(10,25){\line(1,-1){15}}\put(25,10){\line(0,1){15}}\put(25,25){\line(1,-1){15}}\put(40,10){\line(1,1){15}}\put(55,25){\line(0,1){15}}\put(55,40){\line(1,-1){30}}\put(90,30){$\longrightarrow$}\put(110,25){\line(1,-1){15}}\put(125,10){\line(1,1){15}}\put(140,25){\line(0,1){15}}\put(140,40){\line(1,-1){15}}\put(155,25){\line(0,1){15}}\put(155,40){\line(1,-1){30}}\put(187,30){$\longrightarrow$}\put(200,25){\line(1,-1){15}}\put(215,10){\line(1,1){15}}\put(230,25){\line(0,1){15}}\put(230,40){\line(1,-1){30}}\put(260,10){\line(0,1){15}}\put(260,25){\line(1,-1){15}}
	\end{picture}\caption{We aim to shift the leftmost NSE peak to the right twice; however, after the initial move, we encounter a sequence of adjacent peaks. Subsequently, we proceed to move the rightmost peak in this sequence.} \label{figure4}
\end{figure}

If we come up against a sequence of adjacent NSE peaks, then the rightmost peak will become relative height one, and we will move the rightmost NSE peak in the sequence to the right. Increasing the weight of peak $(x,y)$ makes it adjacent to the peak to its right, also commutes the relative height of the two peaks. Figure \ref{figure4} illustrates this situation. Then we can also obtain the following proposition, which is similar to a result of Corteel and Mallet \cite[Proposition 6.5]{cor07}, we omit the proof here.

\begin{pro}	
	The operations in Figure \ref{figure3} and \ref{figure3+}  preserve the number of peaks of relative height $r$ for all $r$.	
\end{pro}

\begin{figure}[h]
	\centering\begin{picture}(340,80)
		\setlength{\unitlength}{0.45mm}\thicklines	\put(10,10){\line(1,1){15}}\put(25,25){\line(1,-1){15}}\put(40,10){\line(1,1){15}}\put(55,25){\line(0,1){15}}\put(55,40){\line(1,-1){30}}\put(90,30){$\longrightarrow$}\put(110,10){\line(1,1){30}}\put(125,25){\line(0,1){2}}\put(140,40){\line(1,-1){15}}\put(155,25){\line(0,1){15}}\put(155,40){\line(1,-1){30}}\put(187,30){$\longrightarrow$}\put(200,10){\line(1,1){30}}\put(215,25){\line(0,1){2}}\put(230,40){\line(1,-1){30}}\put(260,10){\line(0,1){15}}\put(260,25){\line(1,-1){15}}
	\end{picture}\caption{We want to move the leftmost NESE peak to the right twice, but after the first move, we come up against a sequence of
		adjacent peaks. We then move the rightmost peak in this sequence.} \label{figure4+}
\end{figure}
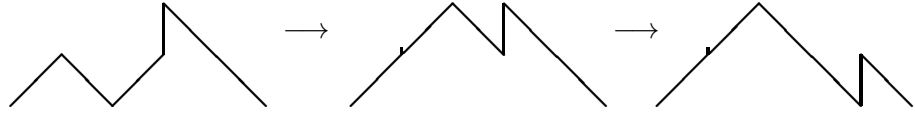
\begin{pro}
		The operations in Figure \ref{figure3} and \ref{figure3+}  do not change the relative position of the adjacent two NSE peaks and NESE peaks.
\end{pro}
The operations in Figure \ref{figure4+} demonstrate this. So,  the NESE peak will stay as the leftmost peak.

We also require a volcanic uplift on the lattice path, which is composed of four types of unit steps.  
The volcanic uplift with only three types of unit steps in the lattice can be seen in \cite{bre87}.

\noindent{\bf The volcanic uplift.}

\begin{figure}[h]
	\centering\begin{picture}(400,80)
		\thicklines 	\setlength{\unitlength}{0.45mm}	\put(10,25){\line(1,-1){15}}\put(25,10){\line(0,1){15}}\put(25,25){\line(1,-1){15}}\put(40,10){\line(1,0){15}}\put(60,25){$\longrightarrow$}\put(85,25){\line(1,-1){15}}\put(100,10){\line(1,1){15}}\put(115,25){\line(0,1){15}}\put(115,40){\line(1,-1){30}}\put(145,10){\line(1,0){15}}
		
		\put(180,10){\line(1,1){15}}\put(195,25){\line(1,-1){15}}\put(210,10){\line(1,0){15}}\put(230,25){$\longrightarrow$}\put(255,10){\line(1,1){30}}\put(270,25){\line(0,1){2}}\put(285,40){\line(1,-1){30}}\put(300,25){\line(0,1){2}}\put(315,10){\line(1,0){15}}
	\end{picture}\caption{The volcanic lift} \label{figure5}
\end{figure}
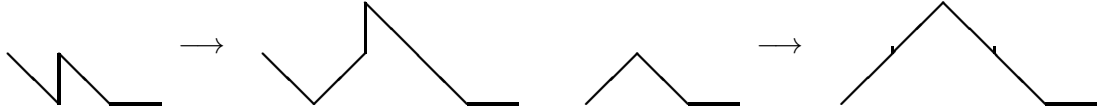

From left to right, we carry out the following steps on each peak one by one. For each NSE peak $(x, y)$, firstly, we delete the two steps adjacent to it and  increase the  $x$-coordinate by 2 for all the steps to the right of $(x+1, 0)$. Then, we place an SE step at $(x, y-1)$, an S step at $(x+1, y)$, an NE step at $(x+1, y+1)$, and an NE step at $(x+2, y)$. 
 Similarly, for an NESE peak,  we delete the two steps adjacent to it and  increase the  $x$-coordinate by 2 for all the steps to the right of $(x+1, 0)$. Then, we place an SE step at $(x-1, y-1)$, an SE step at $(x, y)$, an NE step at $(x+1, y+1)$, and an NE step at $(x+2, y)$. The volcanic uplift is displayed in Figure \ref{figure5}. 

The volcanic uplift successively increases the weight of the peaks by $1,3,5,\ldots, 2s-1$, where $s$ is the number of peaks in the lattice path. Consequently, the major index is increased by $s^2$.

\section{The bijection between overpartition and lattice paths with four kinds of unitary steps}

In the previous section, we referenced some bijections in \cite{chen13a,chen13b,sang15} on overpartitions and defined some operations on lattice paths with four kinds of unitary steps. In this section, we will construct a bijection between overpartitions and lattice paths with four kinds of unitary steps. First, we decompose the Gordon marking of an overpartition $\lambda$ in $\mathcal{U}_{N_1,N_2,\ldots,N_{k-1};a}$, more precisely, $\lambda$ corresponds  to a distinct partition with parts less than $N_1-1$ and $k-1$ partitions with $n_i$ nonnegative parts for $i=1,\ldots,k-1$.

\noindent{\bf The decompsition  of an overpartition Gordon marking} 

Let $\lambda=\lambda_0\in \mathcal{U}_{N_1,N_2,\ldots,N_{k-1};a}$, $\beta=\emptyset$ and $\mu^{(r)}=\emptyset$, for $r=1,\ldots, k-1$. 
\begin{itemize}
\item[Step 1] With the aid of the bijection $\varphi$ introduced in \cite[Section 4]{chen13a}, $\lambda$ is mapped to a distinct partition $\beta$ with parts less than $N_1$ and an overpartition $\alpha_0\in \mathcal{P}_{N_1,N_2,\ldots,N_{k-1};a}$. Let $\alpha=\alpha_0$, for $r$ from $k-1$ to $1$, we proceed to Step 2.
\item[Step 2]
$\alpha$ is an overpartition in $\mathcal{P}_{N_1,N_2,\ldots,N_{r};a}$. Then we apply the map $\psi$ to $\alpha$ and derive an overpartition $\gamma$ in $\mathcal{Q}_{N_1,N_2,\ldots,N_{r};a}$ and a partition $\mu^{(r)}$ in $R_{n_{r}}$. By applying the bijection $\chi$ $n_r$ times, we obtain an overpartition in $\mathcal{P}_{N_1,N_2,\ldots,N_{r-1};a}$, which we denote as $\alpha$. If $r \geq 2$ and $\alpha \neq \emptyset$, we set $r = r - 1$ and return to Step 2; otherwise, we proceed to Step 3.

\item[Step 3] 
We obtain a distinct partition $\beta$ with parts less than $N_1$, and partitions $\mu^{(r)}\in R_{n_r}$, for $r=1,\ldots, k-1$.

\end{itemize}		
\begin{exam}\label{exampar}Given an overpartition $\lambda\in \mathcal{P}_{1,1,5;1}$, where $k=4,a=1$ and 
\[  \lambda=(16,13,12,12,11,\overline{10},\overline{8},8,8,7,\overline{6},6,5,5,4,2,2,\overline{1}),\] then $|\lambda|=136$.
The
Gordon marking of $\lambda$ is
\begin{equation}\label{lambda}
	\lambda=\setcounter{MaxMatrixCols}{16}\begin{bmatrix}
		\ &2&\ &\ &5&\ &7&8&\ &\ &\ &12&\ &\ &\ &\\[3pt]
		\ &2&\ &\ &5&6&\ &8&\ &\ &11&\ &13&\ &\ &\\[3pt]
		\overline{1}&\ &\ &4&\ &\overline{6}&\ &\overline{8}&\ &\overline{10}&\ &12&\ &\ &\ &16
	\end{bmatrix}\begin{matrix}
		3\\[3pt]2\\[3pt]1
	\end{matrix}\;.
\end{equation}
After the first step, we obtain the overpartition $\alpha$ as follows  
\begin{equation*}\label{alpha}
	\setcounter{MaxMatrixCols}{13}\alpha=\begin{bmatrix}
		\ &2&\ &\ &5&6 &\ &8&\ &\ &\ &12&\\[3pt]
		\ &2&\ &4 &\ &6&\ &8&\ &10 &\ &12 &\\[3pt]
		\overline{1}&\ &\ &\overline{4}&\ &\overline{6}&\overline{7}&\ &\ &\overline{10}&\overline{11} &\ &\overline{13}
	\end{bmatrix}\begin{matrix}
		3\\[3pt]2\\[3pt]1
	\end{matrix}\;,
\end{equation*}
and  $\beta=(6,2,1)$.

Since $r=k-1=3$, apply the map $\psi$  to $\alpha$,
then we obtain an overpartition in $P_{7,6,5;1}(121)$  as given below:
\begin{equation*}
	\setcounter{MaxMatrixCols}{13}\begin{bmatrix}\nonumber
\ &2&\ &4 &5 &\ &7&\ &9&\ &\ &\ &\\[3pt]
\ &2&3 &\ &5 &\ &7&8&\ &\  &\ &12 &\\[3pt]
\overline{1}&\ &\overline{3} &\ &\overline{5} &\overline{6} &\ &\overline{8} &\ &\ &\overline{11} &\ &\overline{13}
	\end{bmatrix}\begin{matrix}
		3\\[3pt]2\\[3pt]1
	\end{matrix}\;,
\end{equation*}
together with a partition  $\mu^{(3)}=(7,3,3,3,0)\in R_5$.

Then we apply the bijection $\chi$ to remove the five  $3$-clusters, in order to obtain a new $\alpha$ whose largest label is $2$, as illustrated below.
\begin{equation}
	\setcounter{MaxMatrixCols}{8}\begin{bmatrix}\nonumber
		\ &\ &\ &\  &\  &\ &\ &\ \\[3pt]
		\ &\ &\  &4&\  &\ &\ &\ \\[3pt]
		\ &\ &\  &\overline{4} &\  &\  &\  &\overline{8}
	\end{bmatrix}\begin{matrix}
		3\\[3pt]2\\[3pt]1
	\end{matrix}\;.
\end{equation}

Then for $r=2$, we apply the map $\psi$ to  the aforementioned overpartition to achieve the following overpartition,
\begin{equation}
	\setcounter{MaxMatrixCols}{8}\begin{bmatrix}\nonumber
		\ &\ &\ &\  &\  &\ &\ &\ \\[3pt]
		\ &2&\  &\ &\  &\ &\ &\ \\[3pt]
		\overline{1} &\ &\  &\  &\  &\  &\  &\overline{8}
	\end{bmatrix}\begin{matrix}
		3\\[3pt]2\\[3pt]1
	\end{matrix}\;,
\end{equation}
and a partition $\mu^{(2)}=(5)\in R_1$.

Subsequently, we proceed by employing the bijective function $\chi$ to eliminate the $2$-cluster, thereby yielding an overpartition with the greatest label being $1$, as shown below.
\begin{equation}
	\setcounter{MaxMatrixCols}{8}\begin{bmatrix}\nonumber
		\ &\ &\ &\  &\  &\ &\ &\ \\[3pt]
		\ &\ &\  &\ &\  &\ &\ &\ \\[3pt]
		\  &\ &\  &\  &\  &\  &\overline{7} &\ 
	\end{bmatrix}\begin{matrix}
		3\\[3pt]2\\[3pt]1
	\end{matrix}\;.
\end{equation}
Then, for $r=1$, we apply the map $\psi$ to the above overpartition to obtain the following result
\begin{equation}
	\setcounter{MaxMatrixCols}{8}\begin{bmatrix}\nonumber
		\ &\ &\ &\  &\  &\ &\ &\ \\[3pt]
		\ &\ &\  &\ &\  &\ &\ &\ \\[3pt]
		\overline{1}  &\ &\  &\  &\  &\  &\  &\ 
	\end{bmatrix}\begin{matrix}
		3\\[3pt]2\\[3pt]1
	\end{matrix}\;,
\end{equation}
and a partition $\mu^{(1)}=(6)$.  \qed
\end{exam}
Let $\mathcal{\overline{E}}^{NSE}_{N_1,N_2,\ldots,N_{k-1};a}$ be the set of lattice paths that satisfy the special $(k,a)$-conditions where the leftmost peak is an NSE peak and the number of peaks with relative height greater than or equal to $r$ is $N_r$. Based on the partitions $\beta$ and $\mu^{(r)}$, we now proceed to construct a lattice path $L \in \mathcal{\overline{E}}^{NSE}_{N_1,N_2,\ldots,N_{k-1};a}$ that has $n_r$ peaks with relative height $r$, for $r=1,\ldots, k-1$, and has a major index of $|\lambda_0|$.

\noindent{\bf The construction of the lattice path $L$.} 

Let $L$ be $\emptyset$ and $r=k-1$.

\begin{itemize}
\item[Step 1] 
We insert $n_r$ NSE peaks at $(0, k-r)$, that is, if $r \geq a$, we insert $n_r$ pairs of SE and N steps; if $r < a$, we insert them at $(0, k-a)$.

\item[Step 2] If $a+1\leq r \leq k-1$, we insert an SE step at $(0,k-r+1)$.

\item[Step 3]The partition $\mu^{(r)}$ is a partition that satisfies $\mu^{(r)}_1 \geq \cdots \geq \mu^{(r)}_{n_r} \geq 0$. Then, for $i$ from 1 to $n_r$, we let the weights of the $n_r$ peaks with a relative height of 1 increase successively by $\mu^{(r)}_1, \ldots, \mu^{(r)}_{n_r}$, starting from the rightmost peak to the leftmost peak.

\item[Step 4] If $r\neq 1$, perform the ``volcanic uplift" and let $r=r-1$, then return  to step 1.  If $r=1$, we obtain a lattice path where all  peaks are NSE peaks with $n_r$  peaks having a relative height of $r$, and 
the major index is $|\alpha|$, then proceed to Step 5.

\item[Step 5] $\beta$ is a distinct partition with parts less than $N_1$. Let $\beta$ be $\beta_1\geq \beta_2\geq \cdots \geq \beta_l\geq 1$. Then, for $t$ from $1$ to $l$, we apply the first dilation operation of type A once  and type B dilation operations $\beta_t-1$ times to the lattice path. Then we obtain the lattice path $L$.
\end{itemize}

Through this bijection, an overpartition $\lambda\in \mathcal{U}_{N_1,N_2,\ldots,N_{k-1};a}$ corresponds to a lattice path $L\in \mathcal{\overline{E}}^{NSE}_{N_1,N_2,\ldots,N_{k-1};a}$, and we can also obtain the following results. 

\begin{itemize}
	\item [1.] The weight of the overpartition $\lambda$ is equal to the major index of $L$.
	\item[2.] The number of $r$-cluster in $\lambda$ equals the number of peaks with relative height $r$ in $L$.
	\item[3.] The number of overlined parts in $\lambda$ equals the number of NSE peaks in $L$. 
\end{itemize}

\begin{exam}\label{examlattice}  Based on  Example \ref{exampar}, we construct the corresponding lattice path.

In the case where $r=k-1=3$, in Step 1, we insert $n_3=5$ NSE peaks at the starting point $(0,1)$.
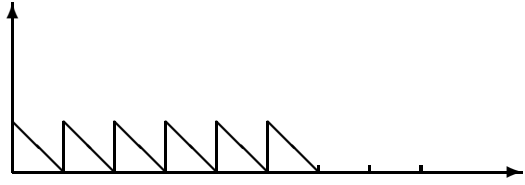
\begin{center}
	\begin{picture}(200,70)
		\thicklines
		\setlength{\unitlength}{0.45mm}
		\multiput(25,10)(15,0){8}{\line(0,1){2}}
		\put(10,10){\vector(1,0){150}}
		\put(10,10){\vector(0,1){50}}
		\multiput(10,25)(15,0){6}{\line(1,-1){15}}
		
		\multiput(25,10)(15,0){5}{\line(0,1){15}}
	\end{picture}
\end{center}
In step 2, since $r=3>a+1=2$, we insert an SE step at $(0,2)$.
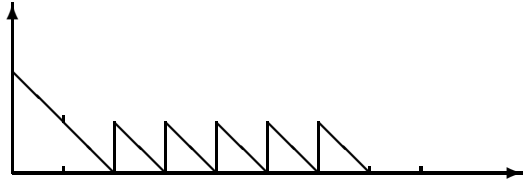
\begin{center}
\begin{picture}(200,70)
\thicklines\setlength{\unitlength}{0.45mm}
		\multiput(25,10)(15,0){8}{\line(0,1){2}}
		\put(10,10){\vector(1,0){150}}
		\put(10,10){\vector(0,1){50}}
		\multiput(25,25)(15,0){6}{\line(1,-1){15}}
		\put(10,40){\line(1,-1){15}}\put(25,25){\line(0,1){2}}
		\multiput(40,10)(15,0){5}{\line(0,1){15}}
	\end{picture}
\end{center}
In step 3, since $\mu^{(3)}=(7,3,3,3,0)$, we increase the weight of each peak successively by $7,3,3,3,0$ starting from the rightmost peak.

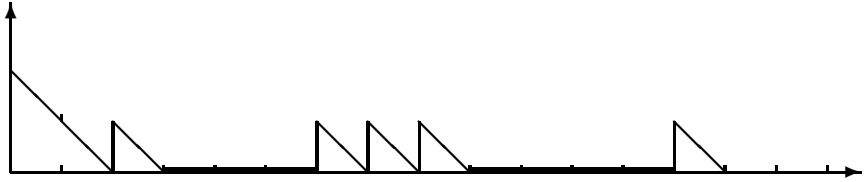
\begin{center}
	\begin{picture}(300,70)
		\thicklines\setlength{\unitlength}{0.45mm}
		\multiput(25,10)(15,0){16}{\line(0,1){2}}
		\put(10,10){\vector(1,0){250}}
		\put(10,10){\vector(0,1){50}}
		\put(10,40){\line(1,-1){15}}\put(25,25){\line(0,1){2}}
	\put(40,10){\line(0,1){15}}\multiput(100,10)(15,0){3}{\line(0,1){15}}\put(205,10){\line(0,1){15}}\multiput(25,25)(15,0){2}{\line(1,-1){15}}\multiput(100,25)(15,0){3}{\line(1,-1){15}}\put(205,25){\line(1,-1){15}}\put(55,11){\line(1,0){45}}\put(145,11){\line(1,0){60}}
	\end{picture}
\end{center}
In Step 4, since $r\neq 1$, perform the volcanic uplift:
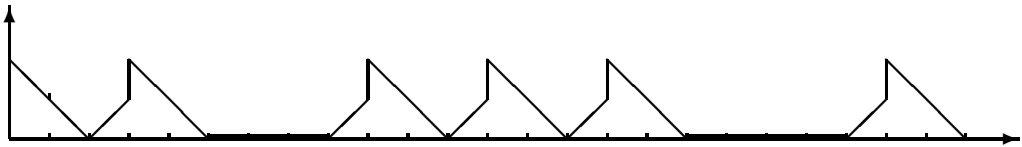
\begin{center}
	\begin{picture}(400,70)
		\thicklines
\multiput(25,10)(15,0){24}{\line(0,1){2}}
\put(10,10){\vector(1,0){380}}
\put(10,10){\vector(0,1){50}}
\put(10,40){\line(1,-1){15}}
\put(25,25){\line(0,1){2}}
\put(40,10){\line(1,1){15}}
\put(55,25){\line(0,1){15}}
\put(55,40){\line(1,-1){30}}
\multiput(145,25)(45,0){3}{\line(0,1){15}}

\put(25,25){\line(1,-1){15}}
\multiput(145,40)(45,0){3}{\line(1,-1){30}}\multiput(130,10)(45,0){3}{\line(1,1){15}}
\put(325,10){\line(1,1){15}}\put(340,25){\line(0,1){15}}\put(340,40){\line(1,-1){30}}
\put(85,11){\line(1,0){45}}
\put(265,11){\line(1,0){60}}
	\end{picture}
\end{center}
Let $r=2$, and go back to Step 1. 
Since $n_2=1$, we add one NSE peak at $(0,2)$.
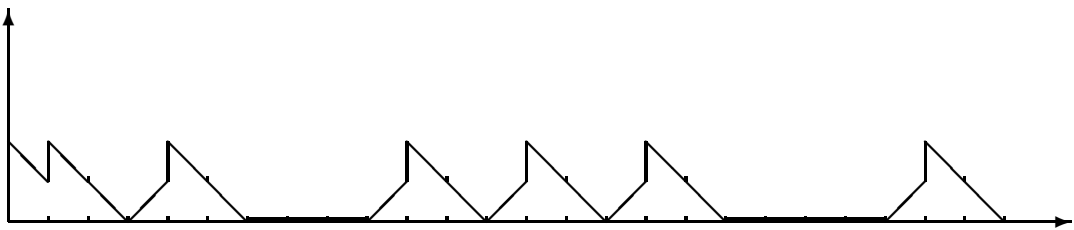
\begin{center}	\begin{picture}(430,90)
		\thicklines
		\multiput(25,10)(15,0){26}{\line(0,1){2}}
		\put(25,10){\vector(1,0){400}}
		\put(25,10){\vector(0,1){80}}
		\put(40,40){\line(1,-1){15}}
		\put(100,25){\line(0,1){2}}
		\put(70,10){\line(1,1){15}}
		\put(85,25){\line(0,1){15}}
		\put(85,40){\line(1,-1){30}}
		\multiput(175,25)(45,0){3}{\line(0,1){15}}
		\multiput(190,25)(45,0){3}{\line(0,1){2}}\put(385,25){\line(0,1){2}}
		\put(55,25){\line(1,-1){15}}
		\multiput(175,40)(45,0){3}{\line(1,-1){30}}\multiput(160,10)(45,0){3}{\line(1,1){15}}
		\put(355,10){\line(1,1){15}}\put(370,25){\line(0,1){15}}\put(370,40){\line(1,-1){30}}
		\put(115,11){\line(1,0){45}}
		\put(295,11){\line(1,0){60}}
		
		\put(25,40){\line(1,-1){15}}\put(40,25){\line(0,1){15}}\put(55,25){\line(0,1){2}}
	\end{picture}
\end{center}
In Step 2, since $r=2\geq a+1$, we add an SE step at $(0,3)$,
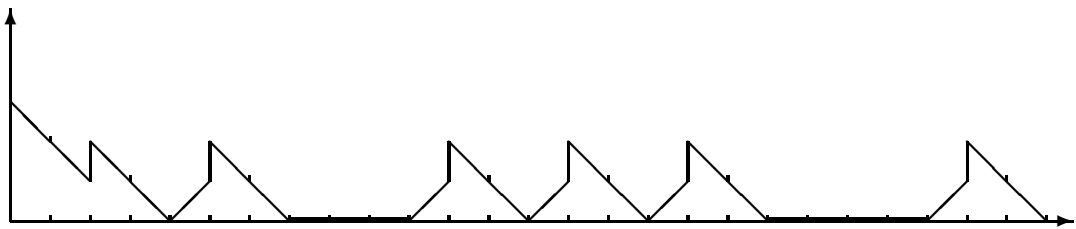
\begin{center}	\begin{picture}(430,90)
		\thicklines
		\multiput(25,10)(15,0){26}{\line(0,1){2}}
		\put(10,10){\vector(1,0){400}}
		\put(10,10){\vector(0,1){80}}
		\put(40,40){\line(1,-1){15}}
		\put(100,25){\line(0,1){2}}
		\put(70,10){\line(1,1){15}}
		\put(85,25){\line(0,1){15}}
		\put(85,40){\line(1,-1){30}}
		\multiput(175,25)(45,0){3}{\line(0,1){15}}
		
		\put(55,25){\line(1,-1){15}}
		\multiput(175,40)(45,0){3}{\line(1,-1){30}}\multiput(160,10)(45,0){3}{\line(1,1){15}}
		\put(355,10){\line(1,1){15}}\put(370,25){\line(0,1){15}}\put(370,40){\line(1,-1){30}}
		\put(115,11){\line(1,0){45}}
		\put(295,11){\line(1,0){60}}
		
		\put(25,40){\line(1,-1){15}}\put(40,25){\line(0,1){15}}\put(55,25){\line(0,1){2}}
		\put(10,55){\line(1,-1){15}}\put(25,40){\line(0,1){2}}
		\multiput(190,25)(45,0){3}{\line(0,1){2}}\put(385,25){\line(0,1){2}}
	\end{picture}
\end{center}
Proceed to  Step 3, $\mu^{(2)}=(5)$, thus we increase the peak with relative height $1$ by $5$, to obtain the following lattice path
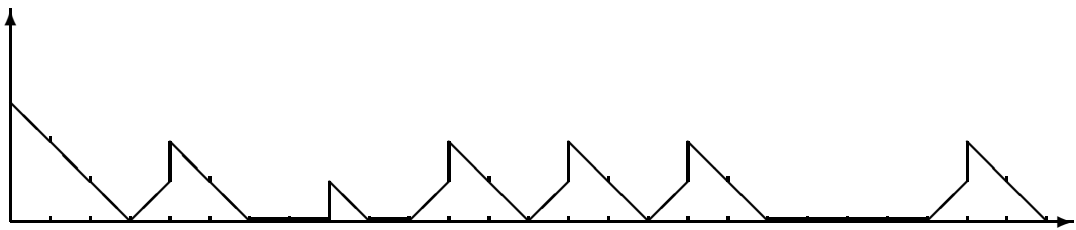
\begin{center}	\begin{picture}(430,90)
		\thicklines
		\multiput(25,10)(15,0){26}{\line(0,1){2}}
		\put(10,10){\vector(1,0){400}}
		\put(10,10){\vector(0,1){80}}
		\put(40,25){\line(1,-1){15}}
		\put(85,25){\line(0,1){2}}
		\put(55,10){\line(1,1){15}}
		\put(130,10){\line(0,1){15}}
		\put(130,25){\line(1,-1){15}}
		\put(85,25){\line(1,-1){15}}
		\multiput(175,25)(45,0){3}{\line(0,1){15}}
		\put(145,11){\line(1,0){15}}
		
		\put(70,40){\line(1,-1){15}}
		\multiput(175,40)(45,0){3}{\line(1,-1){30}}\multiput(160,10)(45,0){3}{\line(1,1){15}}
		\put(355,10){\line(1,1){15}}\put(370,25){\line(0,1){15}}\put(370,40){\line(1,-1){30}}
		\put(100,11){\line(1,0){30}}
		\put(295,11){\line(1,0){60}}
		\multiput(190,25)(45,0){3}{\line(0,1){2}}\put(385,25){\line(0,1){2}}
		\put(25,40){\line(1,-1){15}}\put(70,25){\line(0,1){15}}\put(40,25){\line(0,1){2}}
		\put(10,55){\line(1,-1){15}}\put(25,40){\line(0,1){2}}
	\end{picture}
\end{center}
In Step 4, since $r=2$, perform the volcanic uplift.
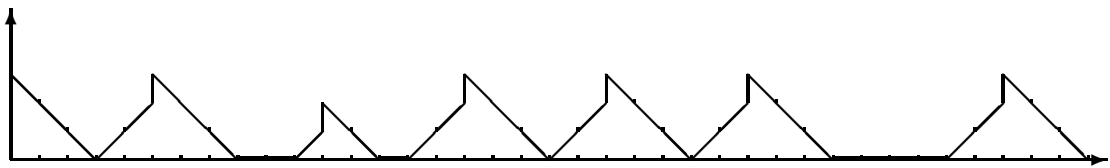
\begin{center}
	\setlength{\unitlength}{0.25mm}	\begin{picture}(650,90)
		\thicklines
		\multiput(25,10)(15,0){38}{\line(0,1){2}}
		\put(10,10){\vector(1,0){580}}
		\put(10,10){\vector(0,1){80}}
		\put(40,25){\line(1,-1){15}}
		\put(115,25){\line(0,1){2}}
		\put(55,10){\line(1,1){30}}
			\put(205,11){\line(1,0){15}}
		\multiput(250,40)(75,0){3}{\line(0,1){15}}
		\put(145,11){\line(1,0){15}}
\put(40,25){\line(0,1){2}}\put(70,25){\line(0,1){2}}		
		\put(85,55){\line(1,-1){45}}
		\multiput(250,55)(75,0){3}{\line(1,-1){45}}
		\put(160,10){\line(1,1){15}}\put(175,25){\line(0,1){15}}\put(175,40){\line(1,-1){30}}\multiput(220,10)(75,0){3}{\line(1,1){30}}
		\put(505,10){\line(1,1){30}}\put(535,40){\line(0,1){15}}\put(535,55){\line(1,-1){45}}
		\put(130,11){\line(1,0){30}}
		\put(445,11){\line(1,0){60}}
		
		\put(25,40){\line(1,-1){15}}\put(85,40){\line(0,1){15}}\put(175,25){\line(0,1){2}}\put(190,25){\line(0,1){2}}\put(235,25){\line(0,1){2}}\put(280,25){\line(0,1){2}}\put(310,25){\line(0,1){2}}\put(355,25){\line(0,1){2}}\put(340,40){\line(0,1){2}}\put(385,25){\line(0,1){2}}\put(430,25){\line(0,1){2}}\put(415,40){\line(0,1){2}}\put(520,25){\line(0,1){2}}\put(565,25){\line(0,1){2}}\put(550,40){\line(0,1){2}}
		\put(10,55){\line(1,-1){15}}\put(25,40){\line(0,1){2}}
	\end{picture}
\end{center}
Set $r=1$, back to Step 1, since $n_1=1$, we add an NSE peak at the beginning.
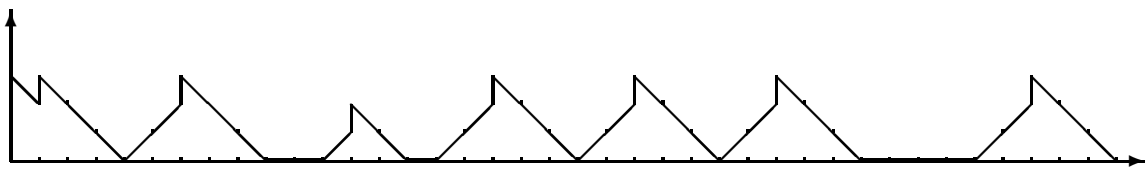
\begin{center}
	\setlength{\unitlength}{0.25mm}	\begin{picture}(650,90)
		\thicklines
		\multiput(25,10)(15,0){39}{\line(0,1){2}}
		\put(10,10){\vector(1,0){600}}
		\put(10,10){\vector(0,1){80}}
		\put(55,25){\line(1,-1){15}}
		\put(85,25){\line(0,1){2}}
		\put(70,10){\line(1,1){30}}
		\put(220,11){\line(1,0){15}}
		\multiput(265,40)(75,0){3}{\line(0,1){15}}
		\put(160,11){\line(1,0){15}}
		
		\put(100,55){\line(1,-1){45}}
		\multiput(265,55)(75,0){3}{\line(1,-1){45}}
		\put(175,10){\line(1,1){15}}\put(190,25){\line(0,1){15}}\put(190,40){\line(1,-1){30}}\multiput(235,10)(75,0){3}{\line(1,1){30}}
		\put(520,10){\line(1,1){30}}\put(550,40){\line(0,1){15}}\put(550,55){\line(1,-1){45}}
		\put(145,11){\line(1,0){30}}
		\put(460,11){\line(1,0){60}}
		
		\put(25,55){\line(1,-1){30}}\put(100,40){\line(0,1){15}}\put(55,25){\line(0,1){2}}
		
		\put(10,55){\line(1,-1){15}}\put(25,40){\line(0,1){15}}
		\put(130,25){\line(0,1){2}}
		\put(205,25){\line(0,1){2}}
		\put(250,25){\line(0,1){2}}
		\put(295,25){\line(0,1){2}}
		\put(280,40){\line(0,1){2}}
		\put(325,25){\line(0,1){2}}
		\put(370,25){\line(0,1){2}}
		\put(355,40){\line(0,1){2}}
		\put(400,25){\line(0,1){2}}
		\put(430,40){\line(0,1){2}}
		\put(445,25){\line(0,1){2}}
		\put(535,25){\line(0,1){2}}
		\put(580,25){\line(0,1){2}}
		\put(565,40){\line(0,1){2}}
		\put(40,40){\line(0,1){2}}
	\end{picture}
\end{center}
In Step 3, since $\mu^{(1)}=(6)$, the leftmost peak with relative height one has its weight increased by $6$.
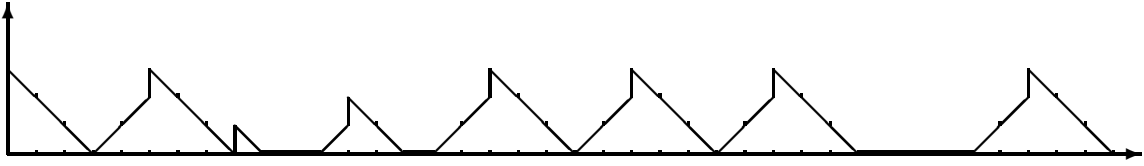
\begin{center}
	\setlength{\unitlength}{0.25mm}	\begin{picture}(650,90)
		\thicklines
		\multiput(25,10)(15,0){39}{\line(0,1){2}}
		\put(10,10){\vector(1,0){600}}
		\put(10,10){\vector(0,1){80}}
		\put(220,11){\line(1,0){15}}
		\put(130,10){\line(0,1){15}}
		\put(55,10){\line(1,1){30}}
		
		\multiput(265,40)(75,0){3}{\line(0,1){15}}
		\put(160,11){\line(1,0){15}}
		
		\put(85,55){\line(1,-1){45}}
		\multiput(265,55)(75,0){3}{\line(1,-1){45}}
		\put(175,10){\line(1,1){15}}\put(190,25){\line(0,1){15}}\put(190,40){\line(1,-1){30}}\multiput(235,10)(75,0){3}{\line(1,1){30}}
		\put(520,10){\line(1,1){30}}\put(550,40){\line(0,1){15}}\put(550,55){\line(1,-1){45}}
		\put(145,11){\line(1,0){30}}
		\put(460,11){\line(1,0){60}}
		\put(10,55){\line(1,-1){45}}\put(85,40){\line(0,1){15}}\put(70,25){\line(0,1){2}}
	\put(130,25){\line(1,-1){15}}
	\put(205,25){\line(0,1){2}}
	\put(250,25){\line(0,1){2}}
	\put(295,25){\line(0,1){2}}
	\put(280,40){\line(0,1){2}}
	\put(325,25){\line(0,1){2}}
	\put(370,25){\line(0,1){2}}
	\put(355,40){\line(0,1){2}}
	\put(400,25){\line(0,1){2}}
	\put(430,40){\line(0,1){2}}
	\put(445,25){\line(0,1){2}}
	\put(535,25){\line(0,1){2}}
	\put(580,25){\line(0,1){2}}
	\put(100,40){\line(0,1){2}}
	\put(115,25){\line(0,1){2}}
	\put(565,40){\line(0,1){2}}
	\put(25,40){\line(0,1){2}}
	\put(40,25){\line(0,1){2}}
	\end{picture}
\end{center}
For $r=1$, there is no volcanic uplift, so proceed to step 5. Recall that $\beta=(6,2,1)$, then by Step 5, we obtain the following lattice path.
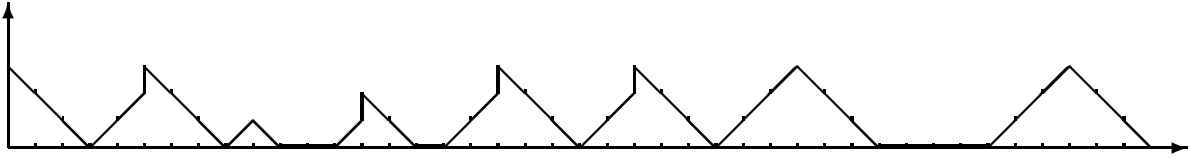
\begin{center}
	\setlength{\unitlength}{0.24mm}	\begin{picture}(660,90)
			\put(25,40){\line(0,1){2}}
		\put(40,25){\line(0,1){2}}
		\thicklines
		\multiput(25,10)(15,0){41}{\line(0,1){2}}
		\put(10,10){\vector(1,0){650}}
		\put(10,10){\vector(0,1){80}}

		\put(55,10){\line(1,1){30}}
		
		\multiput(280,40)(75,0){2}{\line(0,1){15}}
		\put(175,11){\line(1,0){15}}
		\put(85,55){\line(1,-1){45}}
		\multiput(280,55)(75,0){2}{\line(1,-1){45}}
		\put(190,10){\line(1,1){15}}\put(205,25){\line(0,1){15}}\put(205,40){\line(1,-1){30}}
		\put(235,11){\line(1,0){15}}
		\multiput(250,10)(75,0){2}{\line(1,1){30}}\put(400,10){\line(1,1){45}}\put(445,55){\line(1,-1){45}}
		\put(550,10){\line(1,1){45}}\put(595,55){\line(1,-1){45}}
		\put(160,11){\line(1,0){30}}
		\put(490,11){\line(1,0){60}}
		\put(10,55){\line(1,-1){45}}\put(85,40){\line(0,1){15}}\put(70,25){\line(0,1){2}}
		
		\put(145,25){\line(1,-1){15}}\put(130,10){\line(1,1){15}}
	
		\put(220,25){\line(0,1){2}}
		\put(265,25){\line(0,1){2}}
		\put(310,25){\line(0,1){2}}
		\put(295,40){\line(0,1){2}}
		\put(340,25){\line(0,1){2}}
		\put(385,25){\line(0,1){2}}
		\put(370,40){\line(0,1){2}}
		\put(415,25){\line(0,1){2}}
		\put(430,40){\line(0,1){2}}
		\put(475,25){\line(0,1){2}}
			\put(460,40){\line(0,1){2}}
		\put(565,25){\line(0,1){2}}
		\put(625,25){\line(0,1){2}}
		\put(100,40){\line(0,1){2}}
		\put(115,25){\line(0,1){2}}
		\put(580,40){\line(0,1){2}}
			\put(610,40){\line(0,1){2}}
	\end{picture}
\end{center}
\end{exam}

Now we have constructed a lattice path $L$ that satisfies the following conditions: 
\begin{itemize}
	\item[1.] 
	The number of peaks in $L$ is  $7$  which is the number of one marked parts in $\lambda$ and satisfies the special $(k,a)$-conditions with $k=4$ and $a=1$.  
	\item[2.]	
The weights assigned to the peaks are $5,\ 9,\ 13,\ 18,\ 23,\ 29,\ 39$, leading to a major index calculated as the sum $5+9+13+18+23+29+39=136$, which is equal to $|\lambda|$.

\item[3.]
$L$ has $4$ NSE peaks, which coincides with the number of overlined parts in $\lambda$. 

\item[4.] The number of peaks with relative height $3$ in $L$ is $5$, which equals the number of $3$-clusters in $L$. The numbers of peaks with relative heights $2$ and $1$ are both $1$, corresponding to the same number of $2$-clusters and $1$-clusters in $\lambda$, respectively.
\end{itemize}
\begin{thm}\label{thmTU}
All lattice paths in $\mathcal{\overline{E}}^{NSE}_{N_1,N_2,\ldots,N_{k-1};a}$ can be constructed by above steps. Let  $\overline{E}^{NSE}_{N_1,N_2,\ldots,N_{k-1};a}(n)$ be the number of lattice paths in  $\mathcal{\overline{E}}^{NSE}_{N_1,N_2,\ldots,N_{k-1};a}$ whose major index is $n$. Then for $1\leq a\leq k$, we have  
\[\overline{E}^{NSE}_{N_1,N_2,\ldots,N_{k-1};a}(n)=U_{N_1,N_2,\ldots,N_{k-1};a}(n).\]
\end{thm}
\begin{thm}\label{thmTS}
	Let $\mathcal{\overline{E}}^{NESE}_{N_1,N_2,\ldots,N_{k-1};a}$ be the set of lattice paths that satisfy the special $(k,a)$-conditions where the leftmost peak being an NESE peak and the number of peaks with relative height greater than or equal to $r$ is $N_r$. Furthermore, let  $\overline{E}^{NESE}_{N_1,N_2,\ldots,N_{k-1};a}(n)$  denote the number of lattice paths in $\mathcal{\overline{E}}^{NESE}_{N_1,N_2,\ldots,N_{k-1};a}$ that have a major index of $n$.

Then, for $2\leq a\leq k$, we have  \[\overline{E}^{NESE}_{N_1,N_2,\ldots,N_{k-1};a}(n)=\overline{E}^{NSE}_{N_1,N_2,\ldots,N_{k-1};a-1}(n),\]
and 
\begin{equation}\label{ST2}\overline{E}^{NESE}_{N_1,N_2,\ldots,N_{k-1};1}(n)=\overline{E}^{NSE}_{N_1,N_2,\ldots,N_{k-1};k-1}(n-N_1-N_2-\cdots-N_{k-1}).\end{equation}
\end{thm}
\pf  For cases where $a\geq 2$, a simple variation in the construction of a lattice path in $\mathcal{\overline{E}}^{NSE}_{N_1,N_2,\ldots,N_{k-1};a-1}$ can yield a lattice path in $\mathcal{\overline{E}}^{NESE}_{N_1,N_2,\ldots,N_{k-1};a}$. 
To be more precise, 
\begin{itemize}
\item[Step 1.]If the first peak in this path is an NESE peak, we choose the NE step connected to it and change it to an N step.
\item[Step 2.] Move the other Steps on its left one step to the right.
\item[Step 3.] Add an SE step at the starting point to connect to the y-axis, as follows:
\end{itemize}
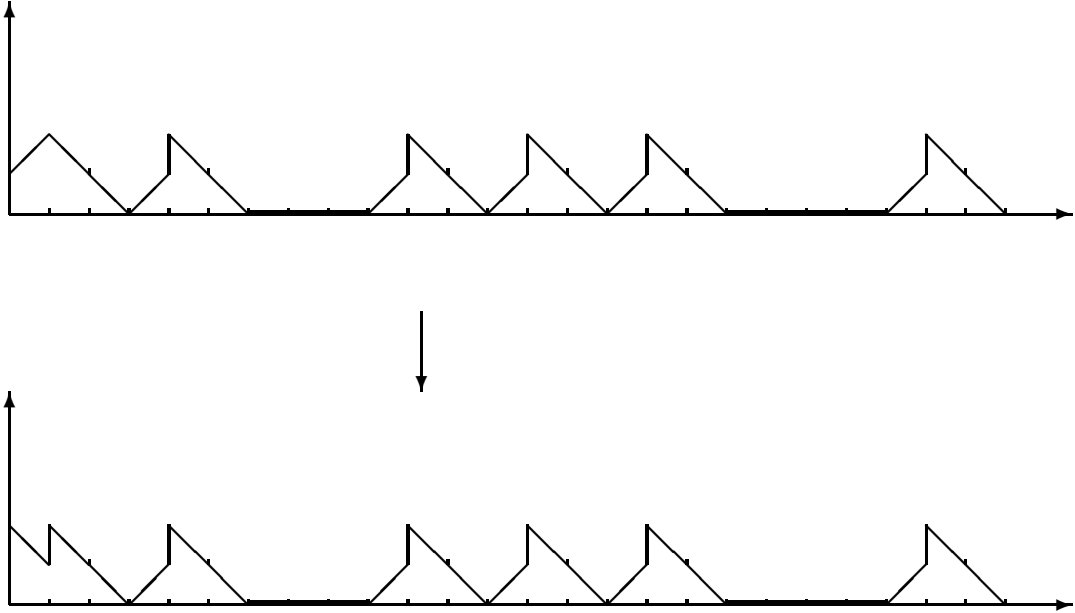
\begin{center}	\begin{picture}(430,90)
		\thicklines
		\multiput(25,10)(15,0){26}{\line(0,1){2}}
		\put(25,10){\vector(1,0){400}}
		\put(25,10){\vector(0,1){80}}
		\put(40,40){\line(1,-1){15}}
		\put(100,25){\line(0,1){2}}
		\put(70,10){\line(1,1){15}}
		\put(85,25){\line(0,1){15}}
		\put(85,40){\line(1,-1){30}}
		\multiput(175,25)(45,0){3}{\line(0,1){15}}
		
		\put(55,25){\line(1,-1){15}}
		\multiput(175,40)(45,0){3}{\line(1,-1){30}}\multiput(160,10)(45,0){3}{\line(1,1){15}}
		\put(355,10){\line(1,1){15}}\put(370,25){\line(0,1){15}}\put(370,40){\line(1,-1){30}}
		\put(115,11){\line(1,0){45}}
		\put(295,11){\line(1,0){60}}
		
		\put(25,25){\line(1,1){15}}\put(55,25){\line(0,1){2}}
	\put(25,40){\line(0,1){2}}
		\multiput(190,25)(45,0){3}{\line(0,1){2}}\put(385,25){\line(0,1){2}}
	\end{picture}
\end{center}
\begin{center}	\begin{picture}(430,130)
		\thicklines
		\put(180,120){\vector(0,-1){30}}
		\multiput(25,10)(15,0){26}{\line(0,1){2}}
		\put(25,10){\vector(1,0){400}}
		\put(25,10){\vector(0,1){80}}
		\put(40,40){\line(1,-1){15}}
		\put(100,25){\line(0,1){2}}
		\put(70,10){\line(1,1){15}}
		\put(85,25){\line(0,1){15}}
		\put(85,40){\line(1,-1){30}}
		\multiput(175,25)(45,0){3}{\line(0,1){15}}
		
		\put(55,25){\line(1,-1){15}}
		\multiput(175,40)(45,0){3}{\line(1,-1){30}}\multiput(160,10)(45,0){3}{\line(1,1){15}}
		\put(355,10){\line(1,1){15}}\put(370,25){\line(0,1){15}}\put(370,40){\line(1,-1){30}}
		\put(115,11){\line(1,0){45}}
		\put(295,11){\line(1,0){60}}
		
		\put(25,40){\line(1,-1){15}}\put(40,25){\line(0,1){15}}\put(55,25){\line(0,1){2}}
		\put(25,40){\line(0,1){2}}
		\multiput(190,25)(45,0){3}{\line(0,1){2}}\put(385,25){\line(0,1){2}}
	\end{picture}
\end{center}
Conversly, for an lattice in $\mathcal{\overline{E}}^{NSE}_{N_1,N_2,\ldots,N_{k};a-1}$, we choose  the left most peak which is an NSE peak, then change the N step  to its left to an NE step 
and delete the first SE step.

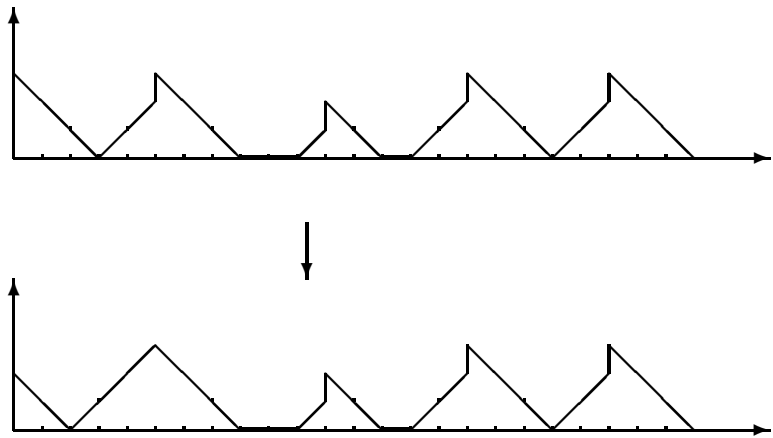
\begin{center}
	\setlength{\unitlength}{0.25mm}	\begin{picture}(350,90)
		\thicklines
		\multiput(25,10)(15,0){24}{\line(0,1){2}}
		\put(25,10){\vector(1,0){400}}
		\put(25,10){\vector(0,1){80}}
		\put(55,25){\line(1,-1){15}}
		\put(85,25){\line(0,1){2}}
		\put(130,25){\line(0,1){2}}
		\put(70,10){\line(1,1){30}}
		\put(220,11){\line(1,0){15}}
		\multiput(265,40)(75,0){2}{\line(0,1){15}}
		\put(160,11){\line(1,0){15}}
		\multiput(205,25)(45,0){3}{\line(0,1){2}}
		\put(100,55){\line(1,-1){45}}
		\multiput(265,55)(75,0){2}{\line(1,-1){45}}
		\put(175,10){\line(1,1){15}}\put(190,25){\line(0,1){15}}\put(190,40){\line(1,-1){30}}\multiput(235,10)(75,0){2}{\line(1,1){30}}
		
		\put(145,11){\line(1,0){30}}

		\put(25,55){\line(1,-1){30}}\put(100,40){\line(0,1){15}}\put(55,25){\line(0,1){2}}
	\end{picture}
\end{center}
\begin{center}
	\setlength{\unitlength}{0.25mm}	\begin{picture}(350,120)
		\thicklines
		\put(180,120){\vector(0,-1){30}}
		\multiput(25,10)(15,0){24}{\line(0,1){2}}
		\put(25,10){\vector(1,0){400}}
		\put(25,10){\vector(0,1){80}}
		\put(55,10){\line(1,1){30}}
		\put(70,25){\line(0,1){2}}
		\put(130,25){\line(0,1){2}}
		
		\put(220,11){\line(1,0){15}}
		\multiput(265,40)(75,0){2}{\line(0,1){15}}
		\put(160,11){\line(1,0){15}}
		\multiput(205,25)(45,0){3}{\line(0,1){2}}
		\put(100,55){\line(1,-1){45}}
		\multiput(265,55)(75,0){2}{\line(1,-1){45}}
		\put(175,10){\line(1,1){15}}\put(190,25){\line(0,1){15}}\put(190,40){\line(1,-1){30}}\multiput(235,10)(75,0){2}{\line(1,1){30}}
		\put(145,11){\line(1,0){30}}
		\put(25,40){\line(1,-1){30}}\put(85,40){\line(1,1){15}}
	\end{picture}
\end{center}

 For the case where $a=1$, we map a lattice path $L\in\mathcal{\overline{E}}^{NSE}_{N_1,N_2,\ldots,N_{k-1},k}$ to a lattice path $L'$ in $\mathcal{\overline{E}}^{NESE}_{N_1,N_2,\ldots,N_{k-1},1}$.  First, in Step 2, for $r$ decreasing from $k-1$ to $2$, one should not do nothing but rather insert an SE step.  Subsequently, let $r$ be the smallest positive integer for which $n_r \geq 1$. In instances similar to when $a \geq 2$, one should insert one NESE peak and $n_r-1$ NSE peaks at $(0, k-r-1)$, as opposed to inserting $n_r$ NSE peaks at $(0,1)$. After completing all the prescribed steps, a lattice path $L' \in \mathcal{\overline{E}}^{NESE}_{N_1,N_2,\ldots,N_{k-1},1}$ can be obtained.
 The insertion step is stated as follows:
 
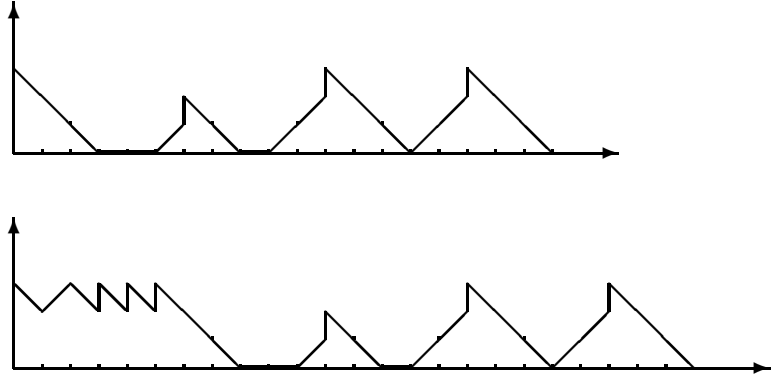
\begin{center}
	\setlength{\unitlength}{0.25mm}	\begin{picture}(500,90)
		\thicklines
		\multiput(100,10)(15,0){20}{\line(0,1){2}}
		\put(100,10){\vector(1,0){320}}
		\put(100,10){\vector(0,1){80}}
		
		\put(130,25){\line(0,1){2}}
		
		\put(220,11){\line(1,0){15}}
		\multiput(265,40)(75,0){2}{\line(0,1){15}}
		\put(160,11){\line(1,0){15}}
		\multiput(205,25)(45,0){3}{\line(0,1){2}}
		\put(100,55){\line(1,-1){45}}
		\multiput(265,55)(75,0){2}{\line(1,-1){45}}
		\put(175,10){\line(1,1){15}}\put(190,25){\line(0,1){15}}\put(190,40){\line(1,-1){30}}\multiput(235,10)(75,0){2}{\line(1,1){30}}		
		\put(145,11){\line(1,0){30}}
		\put(100,40){\line(0,1){15}}\put(55,25){\line(0,1){2}}
	\end{picture}
\end{center}

\begin{center}
	\setlength{\unitlength}{0.25mm}	\begin{picture}(350,90)
		\thicklines
		\multiput(25,10)(15,0){24}{\line(0,1){2}}
		\put(25,10){\vector(1,0){400}}
		\put(25,10){\vector(0,1){80}}
		\put(25,55){\line(1,-1){15}}
	\put(40,40){\line(1,1){15}}	\put(55,55){\line(1,-1){15}}\put(70,40){\line(0,1){15}}\put(70,55){\line(1,-1){15}}\put(85,55){\line(1,-1){15}}\put(85,40){\line(0,1){15}}
		\put(130,25){\line(0,1){2}}
		\put(220,11){\line(1,0){15}}
		\multiput(265,40)(75,0){2}{\line(0,1){15}}
		\put(160,11){\line(1,0){15}}
		\multiput(205,25)(45,0){3}{\line(0,1){2}}
\put(100,55){\line(1,-1){45}}
\multiput(265,55)(75,0){2}{\line(1,-1){45}}
\put(175,10){\line(1,1){15}}\put(190,25){\line(0,1){15}}\put(190,40){\line(1,-1){30}}\multiput(235,10)(75,0){2}{\line(1,1){30}}		
\put(145,11){\line(1,0){30}}
\put(100,40){\line(0,1){15}}
\end{picture}
\end{center}

Now we can see that, in order to construct the lattice path $L'$ in $\mathcal{\overline{E}}^{NESE}_{N_1,N_2,\ldots,N_{k-1},1}$, for $r$ from $k-1$ to $1$ we perform Step 2, which is inserting a total of $k-2$ SE steps and an NESE peak, thereby increasing the major index by  \[n_{k-1}+(n_{k-1}+n_{k-2})+\cdots+(n_{k-1}+n_{k-2}+\cdots+n_1)=N_{k-1}+N_{k-2}+\cdots+N_1,\]
from which we then derive \eqref{ST2}.\qed

\noindent{\bf The proof of Theorem \ref{EB}.} 
Let $m$ denote the number of parts in overpartitions and the sum of relative heights of peaks in lattice paths, thus $m=N_1+N_2+\cdots+N_{k-1}$. Then, by Theorem \ref{thmTS} and \cite[Lemma 3.2]{chen13a}, we obtain the following results.
For the case where $2\leq a\leq k$, we can derive 
\begin{align} \nonumber \overline{E}^{NSE}_{k,a-1}(m,n)&=\sum_{N_1+N_2+\cdots+N_{k-1}=m}\overline{E}^{NSE}_{N_1,N_2,\ldots,N_{k-1};a-1}(n)\\&=\sum_{N_1+N_2+\cdots+N_{k-1}=m}\overline{E}^{NESE}_{N_1,N_2,\ldots,N_{k-1};a}(n)=\overline{E}^{NESE}_{k,a}(m,n)\label{ST}.\end{align}
Similarly, for the case where $a=1$, by employing \eqref{ST2} we  have \begin{equation}\label{FG2}\overline{E}^{NESE}_{k,1}(m,n)=\overline{E}^{NSE}_{k,k-1}(m,n-m).\end{equation}	
Then by Theorem \ref{thmTU}, the result follows  \[\overline{E}^{NSE}_{k,a}(m,n)=U_{k,a}(m,n).\] Moreover, recall the relation derived in \cite[Lemma 3.2]{chen13a}, we reach that  \[\overline{E}^{NESE}_{k,a}(m,n)=V_{k,a}(m,n).\]
Consequently, \[\overline{E}_{k,a}(m,n)=\overline{E}^{NSE}_{k,a}(m,n)+\overline{E}^{NESE}_{k,a}(m,n)=U_{k,a}(m,n)+V_{k,a}(m,n)=\overline{B}_{k,a}(m,n).\]
This completes the proof of Theorem \ref{EB}.  \qed

\section{The proof of Theorem \ref{FD}}

In this section, we shall give the proof of Theorem \ref{FD}.
In the previous section, we have proved Theorem \ref{EB} by constructing a mapping between lattice paths enumerated by $\mathcal{\overline{E}}^{NSE}_{N_1,N_2,\ldots,N_{k-1},a}$ and overpartitions in $\mathcal{U}_{N_1, N_2, \ldots, N_{k-1}; a}(n)$.
Theorem \ref{FD}  can also be derived through a bijection. However, in our proof of Theorem \ref{FD}, we will directly construct the lattice paths enumerated by $\overline{F}_{N_1,N_2,\ldots,N_{k-1};a}(n)$ (the number of lattice paths enumerated by $\overline{F}_{k,a}(n)$ with $n_i$ peaks of relative height $i$), and then we give the generating function for $\overline{F}_{k,a}(n)$. The equality of the generating functions for $\overline{F}_{k,a}(n)$ and $\overline{D}_{k,a}(n)$ will lead to the result in Theorem \ref{FD}.

We first provide the steps for constructing the lattice paths enumerated by $\overline{F}_{N_1,N_2,\ldots,N_{k-1};a}(n)$ in which the leftmost peak is an NSE peak. We denote the set of these lattice paths as $\mathcal{\overline{F}}^{NSE}_{N_1,N_2,\ldots,N_{k-1};a}(n)$.

\noindent{\bf The construction of a lattice path $L\in\mathcal{\overline{F}}^{NSE}_{N_1,N_2,\ldots,N_{k-1};a}(n)$. } 

 For $r$ from $k-1$ to $1$.
\begin{itemize}
\item[Step 1] Insert $n_r=N_r-N_{r+1}$ NSE peaks at the start point  $(0,k-r)$ if $r\geq a$, or at $(0,k-a)$, otherwise.

\item[Step 2] If $a+1\leq r \leq k-1$, an SE step should be introduced at coordinate $(0,k-r+1)$.

\item[Step 3] 
Let the partition $\mu^{(k-1)}$ be any partition with even nonnegative parts such that $\mu^{(k-1)}_1\geq \cdots \geq\mu^{(k-1)}_{n_{k-1}}\geq 0$. For $r=1,\ldots, k-2$, let $\mu^{(r)}$ be any partition with $n_{r}$ nonnegative parts satisfying $\mu^{(r)}_1\geq \cdots \geq\mu^{(r)}_{n_r}\geq 0$. Then, for $i$ from $1$ to $n_r$, we allow the $n_r$ peaks with a relative height of one to successively increase the weights $\mu^{(r)}_1, \ldots ,\mu^{(r)}_{n_r}$ from the rightmost peak to the leftmost peak.

\item[Step 4] 
If $r \neq 1$, execute the volcanic uplift and set $r=r-1$, then return to Step 1. If $r=1$, we obtain a lattice path where all peaks are NSE peaks with $n_r$ of them of relative height $r$  and the major index is $\frac{(N_1+1)N_1}{2}+N_2^2+\cdots+N_{k-1}^2+N_{a+1}+\cdots+N_{k-1}+\sum_{r=1}^{k-1}|\mu^{(r)}|$, proceed to Step 5.

\item[Step 5]
Let $\beta$ be a distinct partition with parts less than $N_1$, with its parts satisfying $N_1>\beta_1>\beta_2>\cdots > \beta_{l}\geq 1$. Then, for $t$ from $1$ to $l$, we apply the first dilation operation of lattice path of Type A once and Type B $\beta_t-1$ times. After these operations, we obtain the lattice path $L$ whose major index is $\frac{(N_1+1)N_1}{2} +N_2^2 +\cdots+ N_{k-1}^2+ N_{a+1} +\cdots+ N_{k-1}+ \sum_{r=1}^{k-1}|\mu^{(r)}| +|\beta|$.

\end{itemize}

\noindent{\bf The construction of a lattice path $L\in\mathcal{\overline{F}}^{NESE}_{N_1,N_2,\ldots,N_{k-1};a}(n)$. } 
\begin{itemize}
\item[Case 1.] For the case $a\geq 2$, 
using the proof method of Theorem \ref{thmTS},  the lattice paths in $\mathcal{\overline{F}}^{NESE}_{N_1,N_2,\ldots,N_{k-1};a}(n)$ can be generated from the lattice paths in $\mathcal{\overline{F}}^{NSE}_{N_1,N_2,\ldots,N_{k-1};a-1}(n)$.

\item[Case 2.] For the case $a=1$, the lattice paths in $\mathcal{\overline{F}}^{NESE}_{N_1,N_2,\ldots,N_{k-1};1}(n)$ can be generated from the lattice paths in $\mathcal{\overline{F}}^{NSE}_{N_1,N_2,\ldots,N_{k-1};k-1}(n)$.

\end{itemize}

We shall prove that $L$ is in $\mathcal{\overline{F}}_{N_1,N_2,\ldots,N_{k-1};a}(n)$. That is to say, we need to verify that $L$ satisfies the condition where the weights $x$ of peaks with the relative height $k-1$ satisfy \eqref{equiv}.

\begin{lem}\label{lemparity}
	All lattice paths in $\mathcal{\overline{F}}_{N_1,N_2,\ldots,N_{k-1};a}(n)=\mathcal{\overline{F}}^{NSE}_{N_1,N_2,\ldots,N_{k-1};a}(n)\bigcup\mathcal{\overline{F}}^{NESE}_{N_1,N_2,\ldots,N_{k-1};a}(n)$ can be generated by the aforementioned steps.
\end{lem}

The peaks with relative height $k-1$ in $L$ are the $n_{k-1}$ peaks that we inserted first. We will first demonstrate that the augmentation of weight in Step 3 does not alter the parity of $x-W(x)$ for these $n_{k-1}$ peaks, where $x$ is the weight of the peaks.

\begin{pro}\label{prodil}
Augmenting the weight of peaks with a relative height of one does not alter the parity of the expression $x-W(x)$ for peaks whose relative height exceeds one, where $x$ denotes the weight of the peak. 
\end{pro}

\begin{figure}[h]\setlength{\unitlength}{0.5mm}\centering\begin{picture}(300,80)\thicklines
\put(9,43){$x_1$}\put(10,10){\line(0,1){25}}\put(10,35){\line(1,-1){25}}\put(35,10){\line(1,1){25}}\put(59,65){$x_2$}\put(60,35){\line(0,1){25}}\put(60,60){\line(1,-1){25}}\put(85,25){$\longrightarrow$}
\put(110,10){\line(1,1){25}}\put(135,35){\line(0,1){25}}\put(134,65){$x_1$}\put(159,65){$x_2$}\put(135,60){\line(1,-1){25}}\put(160,35){\line(0,1){25}}\put(160,60){\line(1,-1){25}}
\put(190,25){$\longrightarrow$}\put(200,10){\line(1,1){25}}\put(225,35){\line(0,1){25}}\put(224,65){$x_1$}\put(274,40){$x_2$}\put(225,60){\line(1,-1){50}}\put(275,10){\line(0,1){25}}\put(275,35){\line(1,-1){25}}
\end{picture}\caption{In the left diagram peak $x_1$ has relative height one and peak $x_2$ has relative height two, in the right diagram their commute the relative height.} \label{figure6}\end{figure}
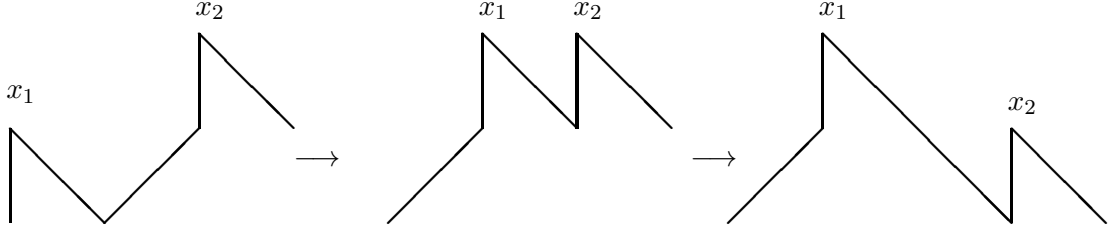

\pf  For $r \leq k-2$, in Step 3, the augmentation of weights is carried out on peaks with relative height one, which will ultimately have a height of $r$ in $L$. When these peaks are inserted in Step 1 and their weights are increased in Step 3 because their relative heights are one, we should note that moving the peaks with relative height one might encounter peaks with a relative height greater than one. As shown in Figure \ref{figure6}, we can see in the left diagram that $x_1$ is the peak with relative height one, and $x_2$ has a relative height of 2. It can be observed that the increase in weight is always performed on the peaks with relative height one. The increase of $x_1$'s weight by two will encounter $x_2$, and then the two peaks will exchange their relative heights. Now we can show that the parity of $x-W(x)$ for the peak with relative height 2 does not change. After the exchange, we can see that the weight of the peak with relative height 2 is decreased by 1, but at the same time, $W(x)$ is also decreased by 1, so the parity of $x-W(x)$ for this peak does not change. Other situations are similar and can be easily checked. We have displayed it in Figure \ref{figure6}.

 \qed

\begin{pro}\label{proins}
	The insertion of an NSE peak does not change the parity of $x-W(x)$ for the peaks to its right.
\end{pro}
\pf
In Example \ref{examlattice}, it is shown that the insertion of an NSE peak will increase the weight of the peaks by one and will also increase the $W(x)$ of these peaks by one, thus the parity of $x-W(x)$ does not change.\qed

\begin{pro}
The first dilation operation of Type A and Type B for
lattice paths does not change the parity of $x-W(x)$ for peaks.
\end{pro}
\pf
The first dilation operation of Type A simply changes the rightmost NSE peak to an NESE peak. 
Then the weight $x$ of the rightmost peak increases by one. Concurrently, the number of NSE peaks with a weight less than or equal to $x$ decreases by one. Therefore, we can observe that the Type A dilation operation does not alter the parity of $x-W(x)$ for all peaks.

The Type B first dilation operation only affects two adjacent peaks, with the left peak changing from an NSE peak to an NESE peak, and the right peak changing from an NESE peak to an NSE peak. Then the weight $x$ of the left peak increases by one and $W(x)$ decreases by one, hence the parity of $x-W(x)$ remains unchanged. The weight $x$ and $W(x)$ of the right peak remain the same. Thus, the parity of $x-W(x)$ is still the same.
\qed
	
\begin{pro}
The volcanic uplift alters the parity of $x-W(x)$ for all peaks.
\end{pro}
\pf The volcanic uplift modifies the parity of the weight associated with each peak, while $W(x)$ remains unchanged. Therefore, after a single volcanic uplift, the parity of $x-W(x)$ for all peaks is altered.\qed

\noindent{\bf The proof of Lemma \ref{lemparity}.}
To prove Lemma \ref{lemparity}, we only need to check whether the peaks with a relative height of $k-1$ satisfy the condition \eqref{equiv}. We first consider the lattice paths in $L\in\mathcal{\overline{F}}^{NSE}_{N_1,N_2,\ldots,N_{k-1};a}(n)$.
The peaks with relative height $k-1$ in $L$ are the $n_{k-1}$ peaks that we inserted first. When these peaks are inserted in Step 1, one can easily check that the parity of $x-W(x)$ is even. And according to Proposition \ref{proins}, the first dilation operation of Type A only changes the rightmost NSE peak to an NESE peak. Then the weight $x$ of the leftmost peak increases by one. At the same time, the number of NSE peaks with a weight less than or equal to $x$ decreases by one. Hence, we can see that the Type A first dilation operation does not change the parity of $x-W(x)$ for all peaks.

The leftmost peak of $L$ is an NSE peak, then for $r$ from $k-1$ to $a+1$, in Step 2, $k-a-1$ SE steps are inserted, which will change the parity of $x-W(x)$ of these $n_{k-1}$ peaks $k-a-1$ times.

For $r \leq k-2$, in Step 3, it has been established that the insertion of additional peaks with a relative height of one does not affect the parity of $x-W(x)$ for the peaks whose relative height exceeds one. For $r=k-1$, the partition $\mu^{(k-1)}$ is a partition with even parts. Then the weight $x$ increases by an even number, while $W(x)$ remains the same. So the parity of $x-W(x)$ also remains unchanged.

We have proven that a single volcanic uplift will change the parity of all peaks. Then we can see that, in Step 4, these $n_{k-1}$ peaks with relative height $k-1$ will undergo $k-2$ volcanic uplifts. So the parity of $x-W(x)$ for these peaks will change $k-2$ times.

In Step 5, we have proven in Proposition \ref{prodil} that the first dilation does not impact the parity of $x-W(x)$ for any peak. 

The leftmost peak of $L$ is an NSE peak, then the parity of $x-W(x)$ for these peaks with relative height $k-1$ has changed $k-a-1$ times in total. Consequently, we obtain that the weights $x$ of the peaks with relative height $k-1$ satisfy 
\[x-W(x)\equiv k-a-1+k-2\equiv a-1\pmod{2},\] 
 which is equivalent to \eqref{equiv}.

Then we consider the case that the lattice path $L'\in \mathcal{\overline{F}}^{NESE}_{N_1,N_2,\ldots,N_{k-1};a}(n)$  can be generated from a lattice path $L\in \mathcal{\overline{F}}^{NSE}_{N_1,N_2,\ldots,N_{k-1};a-1}(n) $ where $a\geq 2$. The peaks in $L$ satisfy that 
\[x-W(x)\equiv a-2\pmod{2},\]  to obtain $L'$, we change the first peak in $L$ to be an $NESE$ peak and the start point change from $(0,k-a+1)$ to $(0,k-a)$. Then $W(x)$ minus $1$. So, in $L'$ we have 
\[x-(W(x)-1)\equiv a-1\pmod{2}.\] 

If $L'\in \mathcal{\overline{F}}^{NESE}_{N_1,N_2,\ldots,N_{k-1};1}(n)$

$L'$ can be generated from $L\in\mathcal{\overline{F}}^{NSE}_{N_1,N_2,\ldots,N_{k-1};k-1}(n-N_1-N_2-\cdots-N_{k-1})$, where the peak $(x,y)$ having an relative height of $k-1$ statisfy that 
\[x-W(x)\equiv k-1-1\pmod{2}.\] Compared to $L$, the construction of $L'$ involves $k-2$ additional insertions of the SE steps, which will change the parity of $x-W(x)$ $k-2$ times. Then we reach that the weights $x$ of  peaks with relative height $k-1$ satisfy 
\[x-W(x)\equiv k-1-1+k-2\equiv 1-1 \pmod{2},\]
which is equivalent to 
\eqref{equiv}.

Combining this result with the proof of Theorem \ref{EB}, we have completed the proof of Lemma \ref{lemparity}.
\qed

It should be noted that we can also obtain the following relations between $\overline{F}^{NSE}_{N_1,N_2,\ldots,N_{k-1};a}$ and $\overline{F}^{NESE}_{N_1,N_2,\ldots,N_{k-1};a}$, which are analogous to those between $\overline{E}^{NSE}_{N_1,N_2,\ldots,N_{k-1};a}$ and $\overline{E}^{NESE}_{N_1,N_2,\ldots,N_{k-1};a}$.
\begin{equation}\label{re1}|\mathcal{\overline{F}}^{NSE}_{N_1,N_2,\ldots,N_{k-1};k-1}(n-N_1-N_2-\cdots-N_{k-1})|=|\mathcal{\overline{F}}^{NESE}_{N_1,N_2,\ldots,N_{k-1};1}(n)|,\end{equation}
and \begin{equation}\label{re2}|\mathcal{\overline{F}}^{NSE}_{N_1,N_2,\ldots,N_{k-1};a-1}(n)|=|\mathcal{\overline{F}}^{NESE}_{N_1,N_2,\ldots,N_{k-1};a}(n)|.\end{equation}

\noindent{\bf The generating functions of $\overline{F}_{k,a}(n)$}

We just need to enumerate the number of lattice paths in $ \mathcal{\overline{F}}^{NSE}_{N_1,N_2,\ldots,N_{k-1};a}(n)$, and then the number of lattice paths in $ \mathcal{\overline{F}}^{NESE}_{N_1,N_2,\ldots,N_{k-1};a}(n)$ will naturally follow.
 
In the first step, the insertion of $n_r$ NSE peaks  results in $n_r$ peaks with the weights $1,2,\ldots,n_r$, and the weight of each original $N_{r+1}=n_{r+1}+ n_{r+2}+\cdots+n_{k-1}$ peak also increases by $n_r$. The major index increases a number of
\[1+2+\cdots+n_r+n_r(n_r+1+\cdots+n_{k-1})=\frac{(n_{r+1})n_r}{2}+n_rN_{r+1}.\]
 
For $r$ from $k-1$ to $1$, the summation is 
\begin{align*}&\frac{(n_{k-1}+1)n_{k-1}}{2}+\frac{(n_{k-2}+1)n_{k-2}}{2}+n_{k-2}N_{k-1}+\cdots+\frac{(n_1+1)n_1}{2}+n_1N_{2}\\&\frac{(n_{k-1}+n_{k-2}+1)(N_{k-1}+n_{k-2})}{2}+\cdots+\frac{(n_1+1)n_1}{2}+n_1N_{2}\\&=\frac{N_{k-2}(N_{k-2}+1)}{2}+\cdots+\frac{(n_1+1)n_1}{2}+n_1N_{2}\\&=\frac{N_1(N_1+1)}{2}
\end{align*}

In Step 2, if $L\in \mathcal{\overline{F}}^{NSE}_{N_1,N_2,\ldots,N_{k-1};a}(n)$, then the insertion of an SE step for $r$ from $k-1$ to $a+1$ results in an increase in the major index by $n_r+\cdots+n_{k-1}=N_r$. In total, the major index has increased by $N_{a+1}+\cdots+N_{k-1}$.

In Step 3, the weight will increase the weight of the partition $\mu^{(r)}$, then for $r\leq k-2$, it is generated by $\frac{1}{(q)_{n_r}}$, as well as $\frac{1}{(q^2;q^2)_{n_{k-1}}}$ for $r=k-1$. Consequently, the factor in the generating function should be,
\[\frac{1}{(q)_{n_1}\cdots(q)_{n_{k-2}}(q^2;q^2)_{n_{k-1}}}.\]

In Step 4, the volcanic uplift will increase the weight of the $N_r=n_r+ \cdots+ n_{k-1}$ peaks by $1, 3, \ldots, 2N_r-1$ each, thus the major index increases by $N_r^2$ for $r$ from $k-1$ to $2$. Therefore, the summation yields 
\[N_2^2+\cdots+N_{k-1}^2.\]

In the last step, the major index has incresed by the number  of $|\beta|$, where $\beta$ is a distinct patition with parts less than $N_1$. Subsequently, $\beta$ can be generated by $(-q)_{N_1-1}$.

Thus, we can obtain the generating function for lattice paths in $\mathcal{\overline{F}}^{NSE}_{N_1,N_2,\ldots,N_{k-1};a}(n)$, and for $1\leq a\leq k-1$, we have 
\[\sum_{L\in \mathcal{\overline{F}}^{NSE}_{N_1,\ldots,N_{k-1};a}(n)}^{n\geq 0}q^n=\sum_{N_1\geq\ldots \geq N_{k-1}\geq 0}\frac{q^{N_1(N_1+1)/2+N_2^2+\cdots+N_{k-1}^2+N_{a+1}+\cdots+N_{k-1}}(-q)_{N_1-1}}{(q)_{n_1}\cdots(q)_{n_{k-2}}(q^2;q^2)_{n_{k-1}}}.\]

Furthermore,
the generating function for lattice paths in $\mathcal{\overline{F}}^{NESE}_{N_1,N_2,\ldots,N_{k-1};a}$ is
\[\sum_{L\in \mathcal{\overline{F}}^{NESE}_{N_1,\ldots,N_{k-1};a}(n)}^{n\geq 0}q^n=\sum_{N_1\geq \ldots \geq N_{k-1}\geq 0}\frac{q^{N_1(N_1+1)/2+N_2^2+\cdots+N_{k-1}^2+N_a+N_{a+1}+\cdots+N_{k-1}}(-q)_{N_1-1}}{(q)_{n_1}\cdots(q)_{n_{k-2}}(q^2;q^2)_{n_{k-1}}},\]
where $1\leq a\leq k-1$.

Consequently, for $1\leq a\leq k-1$, we have the generating function of $\overline{F}_{k,a}(n)$ is 
\begin{equation}
\sum_{n\geq 0}\overline{F}_{k,a}(n)q^n=\sum_{N_1\geq \cdots\geq N_{k-1}\geq 0}\frac{q^{N_1(N_1+1)/2+N_2^2+\cdots+N_{k-1}^2+N_{a+1}+\cdots+N_{k-1}}(-q)_{N_1-1}(1+q^{N_a})}{(q)_{n_1}\cdots(q)_{n_{k-2}}(q^2;q^2)_{n_{k-1}}}.
\end{equation}
Comparing with \eqref{eqD}, we complete the proof of Theorem \ref{FD}. \qed 

\section{The proof of Theorem \ref{thmH} and \ref{thmJ}}

In this section, we shall give the proofs of Theorem \ref{thmH} and \ref{thmJ}. In the last section, we have established that the insertion of peaks, the first dilation operation and the increase by an even number in the weight of peaks will not change the parity of $x-W(x)$ for peaks, while the volcanic uplift will alter the parity of $x-W(x)$ for peaks. To construct the generating functions of $H_{k,a}(n)$ and $J_{k,a}(n)$ step by step, we first discuss how to construct a peak with the desired parity.

\begin{pro}
	The vocanical uplift do not change the parity of \[x-W(x)-\text{relative height of peak}.\]
\end{pro}
\pf 
We have established that volcanic uplift alters the parity of $x-W(x)$, yet concurrently, the relative height increases by one, therefore the parity of $x-W(x)-relative\ height\ of\ peak$ remains unchanged.

\noindent{\bf The construction  of lattice paths in $\mathcal{\overline{E}}^{NSE}_{N_1,N_2,\ldots,N_{k-1};a}$ with peaks of desired parity.}

These constructions should be carried out according to the parities of $k$ and $a$.

\noindent 
For $r$ form $k-1$ to $1$, we shall construct $n_r$ peak  in $L$ of a relative height $r$ with desired parity.

We already have $n_i$ peaks with relative heght $i-r+1$, for $i=r+1,\ldots, k-1$.
\begin{itemize}
	\item[Step 1]Insert $n_r$ NSE peaks, each with a relative height of one, at the initial point. Subsequently, the expression $x-W(x)-1$ has an odd value. 
	\item[Step 2]
If $r \geq a+1$, insert an SE step, which will alter the parity of all peaks. It can be verified that after the peaks have been inserted, there will be $r-a$ SE steps inserted. It is observable that if $r-a$ is odd, then the parity of $x-W(x)-r$ will change an odd number of times.
	
\item[Step 3] This step requires an evaluation of the parity of both $r$ and $a$.
	\begin{itemize}
\item[Case 1.]$r\leq a$. After these peaks have been inserted, no further SE steps will be introduced. The insertion of additional peaks, the increase in the weight of other peaks, volcanic uplift, and the first dilatation operations will not affect the parity of $x-W(x)-the\ relative\ height \ of \ peak$. Therefore, Steps 1, 2, 4, and 5 will render the parity of $x-W(x)-r$ in $L$ odd. Subsequently,  in order to make the peak with a relative height of $r$ to be an even peak, its weight should be increased by any even number. Likewise, if we desire the $r$-peak to be an odd peak, its weight should be increased by an odd number. 
\item[Case 2] 	$r\geq a+1$ and $r-a$ is even. 
Following the insertion of these peaks, an additional $r-a$ SE steps have been  inserted. Subsequently, the procedures denoted as Steps 1, 2, 4, and 5 will maintain the parity of $x-W(x)-the\ relative\ height \ of \ peak$ as odd in $L$. In order to render the peak with a relative height of $r$ as even peak, we may increase its weight by any even number. Likewise, if we wish to have the $r$-peak as an odd peak, we should increase its weight by an odd number.

\item[Case 3.]$r\geq a+1$ and $r-a$ is odd.  After these peaks have been inserted, an additional $r-a$ SE steps have been inserted which alter the parity of $x-W(x)-the\ relative\ height \ of \ peak$. Then, Steps 1, 2, 4, and 5 will change the parity of $x-W(x)-r$ to even in $L$. Then, in order to make the peak with relative height $r$ even, we can increase its weight by any odd number. Similarly, if we wish the $r$-peak to be an odd peak, we should increase its weight by an even number. 
\end{itemize}

\item[Step 4] If $r\neq 1$, proceed with volcanic uplift. It is observable that in this step, once volcanic uplift occurs, the parity of $x-W(x)$ changes once, and likewise, the relative height of the peak increases by one, therefore the parity of $x-W(x)-relative\ height\ of\ peak$ does not change either. Return to Step 1. If $r=1$, proceed to Step 5.

\item[Step 5] 

Let $\beta$ be $\beta_1\geq \beta_2\geq \cdots \geq \beta_{l}\geq 1$, which can be any distinct partition with parts less than $N_1$. Then for $t$ from $1$ to $l$, we apply the first dilation operation once to the lattice path of Type A, and $\beta_t-1$ times to Type B, after which we obtain a lattice path $L$ in $\mathcal{H}^{NSE}_{N_1,N_2,\ldots,N_{k-1};a}$.
\end{itemize}

We can obtain lattice paths with peaks that meet the desired parity conditions by modifying the increment number in Step 3.

\noindent{\bf The construction  of lattice paths in $\mathcal{\overline{E}}^{NESE}_{N_1,N_2,\ldots,N_{k-1};a}$ with desired parity.}

For the case $2 \leq a \leq k-1$, $L' \in \mathcal{\overline{E}}^{NESE}_{N_1,N_2,\ldots,N_{k-1};a}$ can be generated from the lattice path $L \in \mathcal{\overline{E}}^{NSE}_{N_1,N_2,\ldots,N_{k-1};a-1}$. It can be easily checked that the same peaks in $L$ have opposite parity in $L'$.

For the case $a=1$, we construct the path.
Let $r_0$ be the smallst $r$ such that $n_{r}>0$. 
\begin{itemize}
	\item [Step 1.]	
If $r \geq r_0 + 1$, insert $n_r$ NSE peaks with a relative height of one at the coordinate $(0,k-r-1)$. In the case where $r = r_0$, proceed to insert $n_{r_0} - 1$ NSE peaks at $(0,k-r_0)$, and then insert one NESE peak at $(0,k-r_0)$. If $r$ is less than $r_0$ and $n_r$ equals zero, no peaks are to be inserted. It is noteworthy that the insertion of an NESE peak maintains the parity of the peaks to its right.
\item[Step 2.] If $r \geq 1$, insert an SE step, which will commute the parity of all peaks. We can verify that there will be $r-1$ SE steps inserted after the peaks have been inserted. It can be observed that if $r$ is odd, then the parity of $x-W(x)-r$ will change an even number of times.
\item[Step 3] This step requires to account for the parity of $r$.
After these peaks are inserted, $r-1$ SE steps are inserted, the parity of $x-W(x)-r$ is odd in $L$. If $r$ is odd, to make the peak has an even parity, we can increase the weight by any even number. Similarly, if $r$ is even,  we want the $r$-peak to be an even peak, we should increase the weight by an odd number.

Steps 4 and 5 are the same in the construction of the lattice paths in$\mathcal{\overline{E}}^{NSE}_{N_1,N_2,\ldots,N_{k-1};a}$ with desired parity, which do not change the parity of the peaks.\end{itemize}

Let $\mathcal{H}_{N_1,N_2,\ldots,N_{k-1};a}(n)$ (resp. $\mathcal{J}_{N_1,N_2,\ldots,N_{k-1};a}(n)$) denote the set of lattice paths enumerated by $H_{k,a}(n)$ (resp. $H_{k,a}(n)$) where the number of peaks with relative height  $r$ is $n_r$, and set 
\[\mathcal{H}_{N_1,N_2,\ldots,N_{k-1};a}=\bigcup_{n\geq 0}\mathcal{H}_{N_1,N_2,\ldots,N_{k-1};a}(n).\]
Let $\mathcal{H}^{NSE}_{N_1,N_2,\ldots,N_{k-1};a}$ (resp. $\mathcal{H}^{NESE}_{N_1,N_2,\ldots,N_{k-1};a}$) be the subset of $\mathcal{H}_{N_1,N_2,\ldots,N_{k-1};a}$ in which the lattice paths have the NSE peak (resp. NESE peak) as the leftmost peak and let $\mathcal{J}^{NSE}_{N_1,N_2,\ldots,N_{k-1};a}$ (resp. $\mathcal{J}^{NESE}_{N_1,N_2,\ldots,N_{k-1};a}$) be the subset of $\mathcal{J}_{N_1,N_2,\ldots,N_{k-1};a}$ in which the lattice paths have the NSE peak (resp. NESE peak) as the leftmost peak.

\noindent{\bf The generating function of $H_{k,a}(n)$}

We first construct the generating functions for the lattice paths in $\mathcal{H}^{NSE}_{N_1,N_2,\ldots,N_{k-1};a}$. An SE step is inserted if $r\geq a+1$, and to make all the peaks even, as discussed in Step 3, if $r\geq a+1$ and $r-a$ is odd, we should add an odd number to the weight, or add any even number to the weight otherwise.  Subsequent to these considerations, we elucidate the generating function:
\begin{align}
	&\nonumber\frac{q^{N_1(N_1+1)/2+N_2^2+\cdots+N_{k-1}^2+N_{a+1}+\cdots+N_{k-1}+n_{a+1}+n_{a+3}+\cdots}(-q)_{N_1-1}}{(q^2;q^2)_{n_1}\cdots(q^2;q^2)_{n_{k-2}}(q^2;q^2)_{n_{k-1}}}\\
	=&\frac{q^{N_1(N_1+1)/2+N_2^2+\cdots+N_{k-1}^2+2N_{a+1}+2N_{a+3}+2N_{a+5}+\cdots}(-q)_{N_1-1}}{(q^2;q^2)_{n_1}\cdots(q^2;q^2)_{n_{k-2}}(q^2;q^2)_{n_{k-1}}}.
\end{align}

For $a\geq 2$, one can easy get the generating function of  $\mathcal{H}^{NESE}_{N_1,N_2,\ldots,N_{k-1};a}$ is the generating function of  $J^{NSE}_{N_1,N_2,\ldots,N_{k-1};a-1}$:
\begin{align}
&\nonumber\frac{q^{N_1(N_1+1)/2+N_2^2+\cdots+N_{k-1}^2+N_{a}+\cdots+N_{k-1}+n_1+n_2+\cdots+n_{a-1}+n_{a+1}+n_{a+3}+\cdots}(-q)_{N_1-1}}{(q^2;q^2)_{n_1}\cdots(q^2;q^2)_{n_{k-2}}(q^2;q^2)_{n_{k-1}}}\\
&\label{HNE}=\frac{q^{N_1(N_1+1)/2+N_2^2+\cdots+N_{k-1}^2+N_1+2N_{a+1}+2N_{a+3}+2N_{a+5}+\cdots}(-q)_{N_1-1}}{(q^2;q^2)_{n_1}\cdots(q^2;q^2)_{n_{k-2}}(q^2;q^2)_{n_{k-1}}}.
\end{align}	
For $a=1$, one can easily check that also conincident with identity \eqref{HNE}.

Consequently, for $k-1\geq a\geq 1$, the generating function for $H_{k,a}(n)$ can be expressed as:
\small{\begin{align}
		\sum_{n\geq 0}H_{k,a}(n)q^n\label{Hnequiv}=\sum_{n_1,\ldots,n_{k-1}\geq 0}	\frac{q^{N_1(N_1+1)/2+N_2^2+\cdots+N_{k-1}^2+2N_{a+1}+2N_{a+3}+2N_{a+5}+\cdots}(-q)_{N_1}}{(q^2;q^2)_{n_1}\cdots(q^2;q^2)_{n_{k-2}}(q^2;q^2)_{n_{k-1}}}.
\end{align}}

\noindent{\bf The generating function of $J_{k,a}(n)$.}

 Similar to the discussions in the last part, one can readily obtain the generating function for the lattice path in $\mathcal{J}^{NSE}_{N_1,N_2,\ldots,N_{k-1};a}$ is  
	\begin{align*}
		&\frac{q^{N_1(N_1+1)/2+N_2^2+\cdots+N_{k-1}^2+N_{a+1}+\cdots+N_{k-1}+n_1+n_2+\cdots+n_a+n_{a+2}+n_{a+4}+\cdots}(-q)_{N_1-1}}{(q^2;q^2)_{n_1}\cdots(q^2;q^2)_{n_{k-2}}(q^2;q^2)_{n_{k-1}}}\\
		&=\frac{q^{N_1(N_1+1)/2+N_2^2+\cdots+N_{k-1}^2+N_1+2N_{a+2}+2N_{a+4}+2N_{a+6}+\cdots}(-q)_{N_1-1}}{(q^2;q^2)_{n_1}\cdots(q^2;q^2)_{n_{k-2}}(q^2;q^2)_{n_{k-1}}}.
	\end{align*}
	And the  generating function for the lattice path in $\mathcal{J}^{NESE}_{N_1,N_2,\ldots,N_{k-1};a}$ is  
	\begin{align*}
	\frac{q^{N_1(N_1+1)/2+N_2^2+\cdots+N_{k-1}^2+N_{a}+\cdots+N_{k-1}+n_{a}+n_{a+2}+\cdots}(-q)_{N_1-1}}{(q^2;q^2)_{n_1}\cdots(q^2;q^2)_{n_{k-2}}(q^2;q^2)_{n_{k-1}}}\\
		=\frac{q^{N_1(N_1+1)/2+N_2^2+\cdots+N_{k-1}^2+2N_a+2N_{a+2}+2N_{a+4}+2N_{a+6}+\cdots}(-q)_{N_1-1}}{(q^2;q^2)_{n_1}\cdots(q^2;q^2)_{n_{k-2}}(q^2;q^2)_{n_{k-1}}}.
	\end{align*}
	To summarize, we can derive the generating function of $J_{k,a}(n)$ is \small{\begin{align}
			&\nonumber\sum_{n\geq 0}J_{k,a}(n)q^n\nonumber\\&=\sum_{n_1,\ldots,n_{k-1}\geq 0}	\frac{q^{N_1(N_1+1)/2+N_2^2+\cdots+N_{k-1}^2+2N_{a+2}+2N_{a+4}+2N_{a+6}+\cdots}(q^{N_1}+q^{2N_a})(-q)_{N_1-1}}{(q^2;q^2)_{n_1}\cdots(q^2;q^2)_{n_{k-2}}(q^2;q^2)_{n_{k-1}}}.
	\end{align}}

\section{The proofs of Theorem \ref{thmO} and \ref{thmoverO}}

In this section, we will return to the topic of overpartitions. We have provided the bijection between overpartitions in $\mathcal{\overline{B}}_{k,a}(n)$ and lattice paths in $\mathcal{\overline{E}}_{k,a}(n)$. Moreover, in this bijection, each $r$-cluster $C_r(l)$ in $\lambda$ corresponds to a peak $(x,y)$ with relative height $r$. We shall prove the following results.
\begin{pro}\label{pro6}
The cluster $C_r(l)$ has the same parity as the peak $(x,y)$.
\end{pro}	
\pf 
This proposition can be verified through the construction steps for the overpartitions in $\mathcal{\overline{B}}_{k,a}(n)$. We omit it here.

\noindent{\bf The proof of Theorem \ref{thmO}.} Recall that the parity of a cluster is defined as   the opposite parity of the number of even parts in the cluster minus $V(l)$. Then, we can see that  having the number of even parts in each cluster minus $V(l)$ being even is equivalent to all clusters being odd. According to Proposition \ref{pro6}, we can determine that the generating function is \eqref{eqJ}. We have completed the proof. \qed

\noindent{\bf The proof of Theorem \ref{thmoverO}.} 

We should consider how to construct an overpartition such that the clusters within it having the number of odd parts in it minus $V(l)$ is even. From relation \eqref{clusterparity} we konw that is  $\text{the sum of the parts }-V(l)$ should have an even parity.   If the length $r$ of the cluster is odd, then $\text{the sum of the parts}-r-V(l)$ is odd, which means that  the cluster should have an even parity. If $r$ is even, then $\text{the sum of the parts}- r-V(l)$ is even, which means that the cluster should have an odd parity. Therefore, to construct an overpartition enumerated by $\overline{O}_{k,a}(n)$, we should ensure that the clusters satisfy the following conditions:
A cluster with an odd length $r$ should be an even cluster, and a cluster with an even length $r$ should be an odd cluster.

In order to utilize the results and discussion in last section, we back to lattice paths which are corresponding to these overpartitions via the bijection in Section 3.  By the bijection between lattice paths enumerated by $\overline{E}_{k,a}(n)$ and  overpartition enumerated by $\overline{B}_{k,a}(n)$ and Proposition \ref{pro6}, we give the following results.

\begin{thm}\label{thmSO}
	Let $\overline{S}_{k,a}(n)$ denote set of lattice paths that satisfy the special $(k,a)$-conditions, where the peaks of these paths have the parity opposite to that of $r$, with $r$ being the relative height of the peak. Then we have the following results.
	\[\overline{S}_{k,a}(n)=\overline{O}_{k,a}(n).\]
\end{thm}

Based on the discussion in the previous section, we know that, a peak with an odd relative height $r$ should be an even peak, or if $r$ is even, it should be an odd peak. From the discussion in the last section, we can easily derive the generating function for $\overline{S}_{k,a}(n)$.

\begin{itemize}
	\item[1.] When $a$ is odd, the generating function for lattice paths whose leftmost peak is an NSE peak is as follows.
\begin{equation}\sum_{n_1,\ldots,n_{k-1}\geq 0}	\frac{q^{N_1(N_1+1)/2+N_2^2+\cdots+N_{k-1}^2+N_{a+1}+\cdots+N_{k-1}+n_2+n_4+\cdots+n_{a-1}}(-q)_{N_1-1}}{(q^2;q^2)_{n_1}\cdots(q^2;q^2)_{n_{k-2}}(q^2;q^2)_{n_{k-1}}}.
\end{equation}
The generating function for the lattice paths, enumerated by $\overline{S}_{k,a}$, with the leftmost peak being an NESE peak is as follows.
\begin{equation}\sum_{n_1,\ldots,n_{k-1}\geq 0}	\frac{q^{N_1(N_1+1)/2+N_2^2+\cdots+N_{k-1}^2+N_a+N_{a+1}+\cdots+N_{k-1}+n_1+n_3+\cdots+n_{a-2}}(-q)_{N_1-1}}{(q^2;q^2)_{n_1}\cdots(q^2;q^2)_{n_{k-2}}(q^2;q^2)_{n_{k-1}}}.
\end{equation}
Then the generating function is 
\small{\begin{align}
		&\sum_{n\geq 0}\overline{O}_{k,a}(n)q^n\nonumber\\&\label{overOo}=\sum_{n_1,\ldots,n_{k-1}\geq 0}	\frac{q^{N_1(N_1+1)/2+N_2^2+\cdots+N_{k-1}^2+N_{a+1}+\cdots+N_{k-1}}(q^{N_a+n_1+n_3+\cdots+n_{a-2}}+q^{n_2+n_4+\cdots+n_{a-1}})(-q)_{N_1-1}}{(q^2;q^2)_{n_1}\cdots(q^2;q^2)_{n_{k-2}}(q^2;q^2)_{n_{k-1}}}.
\end{align}}

\item[2.] When  $a$ is even, the generating function for lattice paths with the leftmost peak being an NSE peak is as follows.
\small{\begin{equation}
	\sum_{n_1,\ldots,n_{k-1}\geq 0}	\frac{q^{N_1(N_1+1)/2+N_2^2+\cdots+N_{k-1}^2+N_{a+1}+\cdots+N_{k-1}}q^{n_2+n_4+\cdots+n_{a-2}+n_a+n_{a+1}+n_{a+2}+\cdots+n_{k-1}}(-q)_{N_1-1}}{(q^2;q^2)_{n_1}\cdots(q^2;q^2)_{n_{k-2}}(q^2;q^2)_{n_{k-1}}}.
\end{equation}}
Therefore, the generating function for the lattice paths enumerated by $\overline{S}_{k,a}$, with the leftmost peak being an NESE  peak, is as follows.
\small{\begin{equation}
\sum_{n_1,\ldots,n_{k-1}\geq 0}	\frac{q^{N_1(N_1+1)/2+N_2^2+\cdots+N_{k-1}^2+N_a+N_{a+1}+\cdots+N_{k-1}}q^{n_1+n_3+\cdots+n_{a-1}+n_{a+1}+n_{a+2}+\cdots+n_{k-1}}(-q)_{N_1-1}}{(q^2;q^2)_{n_1}\cdots(q^2;q^2)_{n_{k-2}}(q^2;q^2)_{n_{k-1}}}.
\end{equation}}
\small{\begin{align}
		&\sum_{n\geq 0}\overline{O}_{k,a}(n)q^n\nonumber\\\label{overOe}=&\sum_{n_1,\ldots,n_{k-1}\geq 0}	\frac{q^{\frac{N_1(N_1+1)}{2}+N_2^2+\cdots+N_{k-1}^2+N_a+\cdots+N_{k-1}}(q^{N_{a+1}+n_1+n_3+\cdots+n_{a-1}}+q^{n_2+n_4+\cdots+n_{a-2}})(-q)_{N_1-1}}{(q^2;q^2)_{n_1}\cdots(q^2;q^2)_{n_{k-2}}(q^2;q^2)_{n_{k-1}}}.
\end{align}}
\end{itemize}
In conclusion, if $a$ is odd, we can deduce that the generating function for $\overline{S}_{k,a}(n)$ is given by \eqref{overOo}; if $a$ is even, the generating function for $\overline{S}_{k,a}(n)$ is given by \eqref{overOe}. Then, by Theorem \ref{thmSO}, the proof of Theorem \ref{thmoverO} is completed.

\section{The proof of Theorem \ref{thmGT}}

In this section, we will prove Theorem \ref{thmGT}. We know the function \[\frac{q^{N_1(N_1+1)/2+N_2^2+\cdots+N_{k-1}^2}(-q)_{N_1-1}}{(q^2;q^2)_{n_1}\cdots(q^2;q^2)_{n_{k-2}}(q^2;q^2)_{n_{k-1}}}\] can generate lattice paths that satisfy the special $(k,k)$-conditions, where the leftmost peak is an NSE peak and all peaks are odd, as well as lattice paths that satisfy the special $(k,k)$-conditions with the leftmost peak being an NESE peak and all peaks are even, which implies that  all even parity indices are $0$.

We consider the factor $(-yq)_{n_r}$ which generates a distinct partition $(\lambda_1,\ldots,\lambda_l)$ with at most $n_r$ parts, and the exponent of $y$ is the number of parts $l$.

Subsequently, we can observe that for all $r$,  the  lower even parity indices are $0$. The initial increase in weight is from the rightmost $r$-peak to the $\lambda_1$-th $r$-peak from the right. After this initial change,  it can be verified that the parity of these $\lambda_1$ peaks has been altered. In other words, the lower even $r$-parity index is incremented by one. Since $\lambda$ represents a distinct partition, each part $\lambda_i$ changes the parity index once, which implies that the exponent of $y$ is the lower even $r$-peak parity index. Hence, we can ascertain that function \eqref{equ2} serves as the generating function for $G_{k,k}(l,m,n)$.

Based on the bijection between overpartitions and lattice paths outlined in the previous sections, and the correspondence between several parameters, we can easily derive the relationship \eqref{GT} and consequently obtain that \eqref{equ2} also serves as the generating function for $T_{k,k}(l,m,n)$.
This completes the proof of the theorem.

	\vspace{0.5cm}
	\noindent{\bf Acknowledgments.} The author would like to thank Professor Wang Liuquan for his suggestions during the discussions.


\begin{thebibliography}{99} \small
		
		\bibitem{agabre89}	A.K. Agarwal and  D.M. Bressoud, Lattice paths and multiple basic hypergeometric series, Pacific J. Math. 135 (1989)  209–228.
		
		
		\bibitem{and66} G.E. Andrews, An analytic proof of the Rogers-Ramanujan-Gordon
		identities, Amer. J. Math. 88 (1966) 844--846.		
		
		\bibitem{and74}G.E. Andrews, An analytic generalization of the Rogers-Ramanujan
		identities for odd moduli, Proc. Nat. Acad. Sci. USA 71 (1974) 4082--4085.
		
		
		\bibitem{and76}G.E. Andrews, The Theory of Partitions, Cambridge Mathematical Library, Cambridge University Press, Cambridge, 1998. Reprint of Addison-Wesley Publishing Co., 1976.
		
		\bibitem{And10}
		G.E. Andrews, Parity in partition identities, Ramanujan J. 23 (2010) 45--90.		
		
		\bibitem{andbre84} G.E. Andrews and D.M. Bressoud, Identities in combinatorics III: Further aspects of ordered set sorting, Discrete
		Math. 49 (1984) 222--236.
		
		\bibitem{andbre85}	G.E. Andrews and D.M. Bressoud, On the Burge correspondence between partitions and binary words, in: Number Theory,
		Winnipeg, MB 1983, Rocky Mountain J. Math. 15 (2) (1985) 225--233.
		
		
		\bibitem{bre79} D.M. Bressoud, A generalization of the Rogers-Ramanujan identities for all moduli, J. Combin. Theory Ser. A 27 (1979) 64--68.
		
		\bibitem{Bre80} D.M. Bressoud, Analytic and combinatorial generalizations of the Rogers-Ramanujan identities, Mem. Amer. Math. Soc. 24 (1980) 1--54.
		\bibitem{bre87}D.M. Bressoud, Lattice paths and the Rogers–Ramanujan identities, in: Number Theory, Madras 1987, in: Lecture Notes in Math., vol. 1395, Springer, Berlin, 1989, pp. 140--172.
		
		
		
		\bibitem{chen13a}W.Y.C. Chen, D.D.M. Sang and D.Y.H. Shi, The Rogers--Ramanujan--Gordon Theorem for Overpartitions, Proc. London Math. Soc. 106 (3) (2013) 1371--1393.
		
		\bibitem{chen13b}W.Y.C. Chen, D.D.M. Sang and D.Y.H. Shi, An Overpartition Analogue of Bressoud's Theorem of Rogers--Ramanujan Type, Ramanujan J., 36 (2015)  69–80. 
		
		
		\bibitem{cor04}S. Corteel and J. Lovejoy, Overpartitions, Trans. Amer. Math. Soc. 356 (4) (2004) 1623--1635.
		
		
		\bibitem{cor08}S. Corteel, J. Lovejoy and O. Mallet, An extension to overpartitions of the Rogers-Ramanujan identities for even moduli, J. Number Theory 128 (2008) 1602--1621.
		
\bibitem{cor07}S. Corteel and O. Mallet, Overpartitions, lattice paths,
and Rogers–Ramanujan identities, J. Combin. Theory Ser. A 114 (2007) 1407–-1437.			
		
		\bibitem{gor61}B. Gordon, A combinatorial generalization of the Rogers-Ramanujan identities,
		Amer. J. Math. 83 (1961) 393--399.
		\bibitem{kim13}S. Kim and  A.J. Yee,  Rogers-Ramanujan-Gordon identities, generalized G\"{o}llnitz-Gordon identities, and parity questions, J. Comb. Theory, Ser. A 120(5), (2013) 1038--1056.
		
		\bibitem{kur09}K. Kur\c{s}ung\"{o}z,
		Parity considerations in Andrews--Gordon identities, European J. Combin. 31 (2010) 976--1000.
		
		\bibitem{kur10}K. Kur\c{s}ung\"{o}z, Parity considerations in Andrews-Gordon identities, European J. Combin. 31 (2010) 976--1000.		
		
		\bibitem{lov03}J. Lovejoy, Gordon's theorem for overpartitions, J. Combin. Theory Ser. A 103 (2003) 393--401.
		

\bibitem{sang15}D.D.M. Sang and D.Y.H. Shi, An Andrews–Gordon type identity for overpartitions, Ramanujan J, 37 (2015) 653–-679. 

\bibitem{Shi23}R.X.J. Hao and	D.Y.H. Shi, Lattice path and the Rogers-Ramanujan-Gordon type theorem with parity considerations, Adv. Appl. Math.  165 (2025) 102850.
\end{thebibliography}
\end{document}